\documentclass[11pt]{amsart}%
\usepackage{amssymb}
\usepackage{amsfonts}
\usepackage{amsmath}
\usepackage{graphicx}%
\providecommand{\U}[1]{\protect\rule{.1in}{.1in}}
\newtheorem{theorem}{Theorem}
\theoremstyle{plain}
\newtheorem{acknowledgement}{Acknowledgement}

\newtheorem{corollary}{Corollary}

\newtheorem{definition}{Definition}
\newtheorem{example}{Example}

\newtheorem{lemma}{Lemma}

\newtheorem{proposition}{Proposition}
\newtheorem{remark}{Remark}

\numberwithin{equation}{section}
\begin{document}
\title[Stochastic Navier-Sokes equation]{On unbiased stochastic Navier-Stokes equation}
\author{R. Mikulevicius}
\author{B. L. Rozovskii}
\date{May 25, 2010}

\begin{abstract}
A random perturbation of a deterministic Navier-Stokes equation is considered
in the form of an SPDE with Wick type nonlinearity. The nonlinear term of the
perturbation can be characterized as the highest stochastic order
approximation of the original nonlinear term $u{\nabla}u$. This perturbation
is unbiased in that the expectation of a solution of the perturbed equation
solves the deterministic Navier-Stokes equation. The perturbed equation is
solved in the space of generalized stochastic processes using the
Cameron-Martin version of the Wiener chaos expansion. It is shown that the
generalized solution is a Markov process and scales effectively by Catalan numbers.

\end{abstract}
\maketitle

\section{Introduction}

{In this paper we will consider a deterministic Navier-Stokes equation}%
\footnote{Here and below we assume summation over repeating indices in
products}{ }%
\begin{align}
\partial_{t}\mathbf{u}_{0}\left(  t,x\right)   &  =\partial_{i}\left(
a^{\,ij}\left(  t,x\right)  \partial_{j}\mathbf{u}_{0}\left(  t,x\right)
\right) \label{1}\\
&  -u_{0}^{k}\left(  t,x\right)  \partial_{k}\mathbf{u}_{0}%
(t,x)+\mathbf{\nabla}P_{0}\left(  t,x\right)  +\mathbf{f}\left(  t,x\right)
,\nonumber\\
\mathbf{u}_{0}\left(  0,x\right)   &  =\mathbf{w}(x),\operatorname{div}%
\text{\thinspace}\mathbf{u}_{0}=0,\nonumber
\end{align}
{and its stochastic perturbations: }%
\begin{align}
\partial_{t}\mathbf{v}\left(  t,x\right)   &  =\partial_{i}\left(
a^{\,ij}\left(  t,x\right)  \partial_{j}\mathbf{v}\left(  t,x\right)  \right)
-v^{k}\left(  t,x\right)  \partial_{k}\mathbf{v}(t,x)\label{2}\\
&  +\mathbf{\nabla}P\left(  t,x\right)  +\mathbf{f}\left(  t,x\right)
\nonumber\\
&  +[\sigma^{i}(t,x)\partial_{i}\mathbf{v}\left(  t,x\right)  +\mathbf{g}%
\left(  t,x\right)  -\mathbf{\nabla}\tilde{P}\left(  t,x\right)  ]\dot{W}%
_{t},\nonumber\\
\mathbf{v}\left(  0,x\right)   &  =\mathbf{w}\left(  x\right)
,\operatorname{div}\text{\thinspace}\mathbf{v}=0,\nonumber
\end{align}
{and }%
\begin{align}
\partial_{t}\mathbf{u}\left(  t,x\right)   &  =\partial_{i}\left(
a^{\,ij}\left(  t,x\right)  \partial_{j}\mathbf{u}\left(  t,x\right)  \right)
-u^{k}\left(  t,x\right)  \lozenge\partial_{k}\mathbf{u}(t,x)\label{eq:0}\\
&  +\mathbf{\nabla}P\left(  t,x\right)  +\mathbf{f}\left(  t,x\right)
\nonumber\\
&  [\sigma^{i}(t,x)\partial_{i}\mathbf{u}\left(  t,x\right)  +\mathbf{g}%
\left(  t,x\right)  -\mathbf{\nabla}\tilde{P}\left(  t,x\right)  ]\dot{W}%
_{t},\nonumber\\
\mathbf{u}\left(  0,x\right)   &  =\mathbf{w}\left(  x\right)
,\operatorname{div}\text{\thinspace}\mathbf{u}=0,\nonumber
\end{align}
{where }$0\leq t\leq T,x\in\mathbf{R}^{d},d\geq2,${ and }$W_{t}${ is a
cylindrical Wiener process in a separable Hilbert space }$Y$. {The
coefficients }$a^{ij},\sigma^{i}${ and the functions }$\mathbf{f},\mathbf{g}${
are deterministic, }$\sigma^{i}${ and }$\mathbf{g}${ are }$Y${-valued. }Symbol
$\lozenge$ stands for Wick product (see Section {\ref{WickSkor} and references
\cite{GJ}, \cite{Simon})}. Wick product is a stochastic convolution. It could
be interpreted as a generalized Malliavin divergence operator with respect to
Gaussian measure associated with white noise $\dot{W}$ (see \cite{LRS}).

Stochastic PDEs involving Wick product type nonlinearity were originally
discussed in the literature related to the Parisi-Wu program (see
\cite{albr9}, \cite{Laz} and also {\cite{MAR} (Section 6)}). In these papers
Wick product was defined by (Gaussian) invariant measures for the related
PDEs. Other related papers include: \cite{dpD1}, \cite{kon2}{, \cite{Luo2},
\cite{LR2}-\cite{LRS}, \cite{wrk}, etc.). }

Equations {(\ref{2})} and {(\ref{eq:0}) }are \textit{stochastic perturbations}
of the deterministic Navier-Stokes {equation (\ref{1}). }It is shown in
Section 3 that the generalized{ mean, i.e. the zero-order coefficient in the
Wiener chaos expansion (\ref{f2}) of the solution of equation (\ref{eq:0}), is
a solution of equation (\ref{1}), i.e. }%
\begin{equation}
\mathbf{Eu}\left(  t,x\right)  =\mathbf{u}_{0}\left(  t,x\right)  ,
\label{eq:unbaised}%
\end{equation}
where $\mathbf{u}_{0}\left(  t,x\right)  $ is a solution of {equation
(\ref{1}).} {In other words, the solution of stochastic Navier-Stokes
equation} {(\ref{eq:0}) is a }mean preserving{ (\textit{unbiased}) random
perturbation of deterministic Navier-Stokes equation (\ref{1}).\textit{ }}

{ Obviously, this nice property does not hold for equation (\ref{2}) }or other
standard stochastic perturbations of Navier-Stokes equation (e.g. random
initial conditions, random forcing, etc.)

{In fact, equation (\ref{eq:0}) could be viewed as an \textit{approximation}
of stochastic Navier-Stokes equation (\ref{2}). Indeed, under certain natural
assumptions, the following equality holds: }%
\begin{equation}
\mathbf{v}\nabla\mathbf{v}=\sum_{n=0}^{\infty}\frac{\mathcal{D}^{n}%
\mathbf{v}\lozenge\mathcal{D}^{n}\nabla\mathbf{v}}{n!} \label{eq: MDexpansion}%
\end{equation}
where $\mathcal{D}^{n}$ is the $n^{\text{\textrm{th}}}$ power of Malliavin
derivative $\mathcal{D}$. Taking into account expansion (\ref{eq: MDexpansion}%
), $\mathbf{v}\lozenge\nabla\mathbf{v}$ could be viewed as an approximation of
the product $\mathbf{v}\nabla\mathbf{v}$. In fact, $\mathbf{v}\lozenge
\nabla\mathbf{v}$ is the highest stochastic order approximation of
$\mathbf{v}\nabla\mathbf{v}$ { (see Appendix I, Proposition
\ref{prop:wickvsdot} and Remark \ref{re:hiorder})}.

{Stochastic Navier-Stokes equation (\ref{2}) is reasonably well understood and
there exists substantial literature on its analytical properties as well as
its derivation from the first principles (see e.g. \cite{mr22}, \cite{mr23}
and the references therein). In this paper we will be focusing mostly on
equation (\ref{eq:0}).}

{Burger's equation with Wick product was considered in \cite{kon3},
\cite{holoks}, and \cite{kl}, see also the references therein. }

{It was shown in \cite{mr23} that under reasonable assumptions stochastic
Navier-Stokes equation (\ref{2}) has a square integrable solution. Moreover,
this solution can be formally written in the Wiener chaos expansion form: }%
\[
\mathbf{v}\left(  t,x\right)  =\sum_{\alpha}\mathbf{v}_{\alpha}\left(
t,x\right)  \xi_{\alpha},
\]
{where }$\left\{  \xi_{\alpha},\alpha\in J\right\}  ${ is the Cameron-Martin
basis generated by }$\dot{W}_{t}${, }$\mathbf{v}_{\alpha}\left(  t,x\right)
=\mathbf{E}\left(  \mathbf{v}\left(  t,x\right)  \xi_{\alpha}\right)  ${, and
}$J${ is the set of multiindices }$\alpha=\left\{  \alpha_{k},\text{ }%
k\geq1\right\}  ${ such that for every }$k,${ }$\alpha_{k}\in\mathbf{N}%
_{0}(\mathbf{N}_{0}\mathbf{=}\left\{  0,1,2,\ldots\right\}  )${ and }%
$|\alpha|=\sum_{k}\alpha_{k}<\infty.$

{It was shown in \cite{mr23} that the Wiener chaos coefficients }%
$\mathbf{v}_{\alpha}\left(  t,x\right)  ${\ satisfy the propagator
equation:}\newline%
\begin{align}
\partial_{t}\mathbf{v}_{\alpha}\left(  t,x\right)   &  =\partial_{i}\left(
a^{\,ij}\partial_{j}\mathbf{v}_{\alpha}\left(  t,x\right)  \right)  -\nabla
P\left(  t,x\right)  +\mathbf{f}\left(  t,x\right)  I_{\left\{  \left\vert
\alpha\right\vert =0\right\}  }\label{eq: propv}\\
&  -\sum_{p}\sum_{0\leq\beta\leq\alpha}c(\alpha,\beta,p)\left(  \mathbf{v}%
_{\beta+p},\nabla\right)  \mathbf{v}_{\alpha+p-\beta}\left(  t,x\right)
\nonumber\\
&  \sum_{k}\sqrt{\alpha_{k}}[(\sigma^{i},e_{k})_{Y}\partial_{i}\mathbf{v}%
_{\alpha(k)}\left(  t,x\right)  +I_{\left\{  \left\vert \alpha\right\vert
=1\right\}  }(\mathbf{g},e_{k})_{Y}\mathbf{]};\nonumber\\
\mathbf{v}_{\alpha}(0,x)  &  =\mathbf{w}_{\alpha}(x),\operatorname{div}%
\text{\thinspace}\mathbf{v}_{\alpha}=0,\nonumber
\end{align}
{where }$\alpha(k)=\left(  \alpha_{1},\alpha_{2},...\alpha_{k-1},\alpha
_{k}-1,\alpha_{k+1},...\right)  ${ and \ \ }%
\[
c(\alpha,\beta,p)=\left[  \binom{\alpha}{\beta}\binom{\beta+p}{p}\binom
{\alpha+p-\beta}{p}\right]  ^{1/2}.
\]

{\bigskip One advantage of the Wiener chaos representation is that it provides
convenient explicit formulae for computing statistical moments of the random
field }$\mathbf{u}(t,x)${ via Wiener chaos coefficients (see \cite{mr22},
\cite{mr23}). For example, \ }%
\[%
\begin{array}
[c]{c}%
\mathbf{E}u^{i}(t,x)=u_{0}^{i}(t,x),\;\mathbf{E}\left(  u^{i}(t,x)u^{^{_{j}}%
}(t,y)\right)  =\sum_{|\alpha|<\infty}u_{\alpha}^{i}(t,x)u_{\alpha}^{j}(t,y),
\end{array}
\]

{In this paper, we prove that the WCE coefficients of a solution of the
unbiased stochastic Navier-Stokes equation (\ref{eq:0}) are given by }%
\begin{align}
\partial_{t}\mathbf{u}_{\alpha}\left(  t,x\right)   &  =\partial_{i}\left(
a^{\,ij}\partial_{j}\mathbf{u}_{\alpha}\left(  t,x\right)  \right)  -\nabla
P\left(  t,x\right)  +\mathbf{f}\left(  t,x\right)  I_{\left\{  \left\vert
\alpha\right\vert =0\right\}  }\label{eq: propu}\\
&  -\sum_{0\leq\beta\leq\alpha}\sqrt{\binom{\alpha}{\beta}}\left(
\mathbf{u}_{\alpha-\beta},\nabla\right)  \mathbf{u}_{\beta}\left(  t,x\right)
\nonumber\\
&  +\sum_{k}\sqrt{\alpha_{k}}[\left(  \sigma^{j}(t,x),e_{k}(t)\right)
_{Y}\partial_{j}\mathbf{u}_{\alpha(k)}\left(  t,x\right) \nonumber\\
&  +\left(  \mathbf{g}\left(  t,x\right)  ,e_{k}(t)\right)  _{Y}1_{|\alpha
|=1}]\nonumber\\
\mathbf{u}_{\alpha}(0,x)  &  =\mathbf{w}(x),\operatorname{div}\text{\thinspace
}\mathbf{u}_{\alpha}=0.\nonumber
\end{align}

{Clearly, this system of equations is much simpler than equation
(\ref{eq: propv}). }

{If }$\alpha=0${, then }$\mathbf{u}_{\alpha}(t,x)${ is a solution of
deterministic Navier-Stokes equation (\ref{1}). The remaining components are
governed by Stokes equations and could be solved sequentially. From the
computational point of view this is a substantial advantage. Indeed, the
propagator for equation (\ref{2}) is a full nonlinear system while equation
(\ref{eq:0}) is a lower triangular system and only the first equation of this
system is nonlinear.}

{An important feature of equation (\ref{eq:0}) is that }%
\[
\left(  u^{k}\left(  t\right)  \lozenge\partial_{k}\mathbf{u}(t),\mathbf{u}%
\left(  t\right)  \right)  _{L_{2}(\mathbf{R}^{d})}\neq0.
\]
{Therefore, one could not expect a solution of (\ref{eq:0}) to be square
integrable. This effect is not specific to stochastic Navier-Stokes equation.
In fact, it is common for a large class of stochastic bilinear PDEs (see e.g.
\cite{LR2}, \cite{LR3}).}

{In this paper we }consider equation (\ref{eq:0}) in the class of formal
Wiener chaos expansions and show that a formal series%
\begin{equation}
\mathbf{u}\left(  t,x\right)  =\sum_{\alpha\in J}\mathbf{u}_{\alpha}\left(
t,x\right)  \xi_{\alpha} \label{f2}%
\end{equation}
solves (\ref{eq:0}) if and only if \textbf{$u$}$_{\alpha}\left(  t,x\right)
${ are given by equation (\ref{eq: propu}). To make this solution square
integrable, }we rescale it {using second quantization operators (see Appendix
I, 5.1). }It is shown that \textbf{$u$}$(t,x)$ is the limit of square
integrable solutions of the rescaled equations (see Proposition \ref{p1}).

{Convergence of this solution is determined by a system of positive weights
}$\left\{  r_{\alpha}\right\}  _{\left\vert \alpha\right\vert <\infty}${ such
that }$\ ${ }%
\begin{equation}
\left\Vert \mathbf{u}\right\Vert _{\mathcal{R}}^{2}:=\sum_{\left\vert
\alpha\right\vert <\infty}r_{\alpha}^{2}\left\Vert \mathbf{u}_{\alpha}\left(
t\right)  \right\Vert _{L_{2}\left(  (0,T);\mathbf{R}^{d}\right)  }^{2}%
<\infty. \label{eq: stnorms}%
\end{equation}
It turned out, that {Catalan numbers (see \cite{catalan}, \cite{kl}) are
critical for an appropriate choice of the weights }$r_{\alpha}${ in
(\ref{eq: stnorms}) (see Proposition \ref{cle1}). }

In addition, it was shown that{ a solution of equation (\ref{eq:0}) belongs to
the intersection of Sobolev spaces }$\mathbb{H}_{2}^{2}\left(  \mathbf{R}%
^{d}\right)  \cap\mathbb{H}_{p}^{2}\left(  \mathbf{R}^{d}\right)  ${ for
}$p>d${. We have also demonstrated that uniqueness of a solution of equation
(\ref{eq:0}) holds under the same assumptions that guarantee uniqueness for
the related deterministic Navier-Stokes equation. }Although $\xi_{\alpha}$ in
(\ref{f2}) are not $(\mathcal{F}_{t}^{W})$-adapted, we{ prove that the
generalized solution is }$(\mathcal{F}_{t}^{W})$-adapted {and Markov (see
Theorem \ref{main} and Corollary \ref{MP}) . }

{It is not clear how, if at all, the unbiased Navier-Stokes equation fits into
classical fluid mechanics. Nevertheless, equation (\ref{eq:0}) is
\textquotedblleft physical\textquotedblright\ in that it could be derived from
the second Newton law (under appropriate assumptions on the velocity field),
much the same way as the classical Navier-Stokes equation (see Appendix I,
5.2). Also, it was shown recently (see \cite{kkrv}) that, after Catalan type
rescaling, finite dimensional projections of unbiased Navier-Stokes equation
present an accurate and numerically inexpensive approximation of stochastic
Navier-Stokes equation (\ref{2})}

We conclude this section with an outline of some notations that will be used
in the paper.

\subsection{Notation}

Let us fix a separable Hilbert space $Y$. The scalar product of $x,y\in$%
\ $Y$\ will be denoted $\ $by $(x,y)_{Y}.$

If $u$\ is a function on $\mathbf{R}^{d},$\ the following notational
conventions will be used for its partial derivatives: $\partial_{i}u=\partial
u/\partial x_{i},\partial_{ij}^{2}=\partial^{2}u/\partial x_{i}\partial
x_{j},$\ $\partial_{t}u=\partial u/\partial t$, and $\nabla u=\partial
u=(\partial_{1}u,\ldots,\partial_{d}u),$ and $\partial^{2}u=(\partial_{ij}%
^{2}u)$\ denotes the Hessian matrix of second derivatives. Let $\alpha=\left(
\alpha_{1},...,\alpha_{d}\right)  $\ be a multi-index, $\alpha_{i}%
\in\mathbf{N}_{0}\mathbf{=}\left\{  0,1,2,\ldots\right\}  ,i=1,\ldots,d$, then
$\partial_{x}^{\alpha}=\Pi_{i=1}^{d}\partial_{x_{i}}^{\alpha_{i}}.$

Vector fields on $\mathbf{R}^{d}$ are denoted by boldface letters. This
convention also applies if the entries of the vector field are taking values
in a Hilbert space.

We denote $\mathbf{N}=\left\{  1,2,\ldots\right\}  $.

For a Banach space $E$, we denote $C\left(  [0,T],E\right)  $ the space of
continuous $E$-valued functions.

$C_{0}^{\infty}=C_{0}^{\infty}(\mathbf{R}^{d})$\ denotes the set of all
infinitely differentiable functions on $\mathbf{R}^{d}$\ with compact support.

For $s\in(-\infty,\infty),$ write $\Lambda^{s}=\Lambda_{x}^{s}=\left(
1-\sum_{i=1}^{d}\partial^{2}/\partial x_{i}^{2}\right)  ^{s/2}.$ For
$p\in\lbrack1,\infty)$\ and $s\in(-\infty,\infty)$, we define the space
$H_{p}^{s}=H_{p}^{s}(\mathbf{R}^{d})$\ as the space of generalized real valued
functions $u$\ with the finite norm
\[
|u|_{s,p}=|\Lambda^{s}u|_{p},
\]
where $|\cdot|_{p}$\ is the $L_{p}$\ norm. Obviously, $H_{p}^{0}=L_{p}$. Note
that if $s\geq0$\ is an integer, the space $H_{p}^{s}$\ coincides with the
Sobolev space $W_{p}^{s}=W_{p}^{s}(\mathbf{R}^{d})$.

The spaces $C_{0}^{\infty}(\mathbf{R}^{d}),\,H_{p}^{s}\left(  \mathbf{R}%
^{d}\right)  $ can be extended to vector functions (denoted by bold-faced
letters). For example, the space of all vector functions \ $\mathbf{u}%
=(u^{1},\ldots,u^{d})$\ such that $\Lambda^{s}u^{l}\in L_{p},$\ $l=1,\ldots
,d$, with the finite norm
\[
|\mathbf{u}|_{s,p}=(\sum_{l}|u^{l}|_{s,p}^{p})^{1/p},
\]
is denoted by $\mathbb{H}_{p}^{s}=\mathbb{H}_{p}^{s}(\mathbf{R}^{d})$.
Similarly, we denote by $\mathbb{H}_{p}^{s}(Y)=\mathbb{H}_{p}^{s}%
(\mathbf{R}^{d},Y)$\ the space of all vector functions $\mathbf{g}%
=(g^{l})_{1\leq l\leq d},$\ with $Y$-valued components $g^{l},$\ $1\leq l\leq
d,$ so that $||\mathbf{g}||_{s,p}=(\sum_{l}|g^{l}|_{s,p}^{p})^{1/p}<\infty$.
Also, for brevity, the norm $||\mathbf{g}||_{0,p}$ is denoted by
$||\mathbf{g}||_{p}.$

When $s=0,$\ $\mathbb{H}_{p}^{s}(Y)=\mathbb{L}_{p}(\,Y)=\mathbb{L}%
_{p}(\mathbf{R}^{d},\,Y).$ To forcefully distinguish $L_{p}$-norms in spaces
of $Y$-valued functions, we write $||\cdot||_{p},$\ while in all other cases a
norm is denoted by $\left\vert \cdot\right\vert _{p}.$ The duality
$\left\langle \cdot\mathbf{,\cdot}\right\rangle _{s}$ between \textrm{\ }%
$\mathbb{H}_{q}^{s}\left(  \mathbf{R}^{d}\right)  ,$ and $\mathbb{H}_{p}%
^{-s}\left(  \mathbf{R}^{d}\right)  $ where\textrm{\ }$p\geq2$\textrm{\ }and
$q=p/\left(  p-1\right)  $ is defined by
\[
\left\langle \boldsymbol{\phi,\psi}\right\rangle _{s}=\left\langle
\boldsymbol{\phi,\psi}\right\rangle _{s,p}=\sum_{i=1}^{d}\int_{\mathbf{R}^{d}%
}(\Lambda^{s}\phi^{i})\left(  x\right)  (\Lambda^{-s}\psi^{i})\left(
x\right)  dx,\mathbf{\phi}\in\mathrm{\ }\mathbb{H}_{q}^{s},\mathbf{\psi}%
\in\mathbb{H}_{p}^{-s}.
\]

\section{{Generalized Random Variables and Processes}}

\subsection{{\label{ss1}Wiener Chaos}}

{To begin with, we shall introduce some basic notation and recall a few
fundamental facts of infinite-dimensional stochastic calculus. Let us fix a
separable Hilbert spaces }$Y${ and }$\mathbf{H}=L_{2}([0,T],Y).${ Let }%
$\{\ell_{i},i\geq1\}${ be a complete orthonormal basis (CONS) in }$Y${
and}$\{m_{i},i\geq1\}${ be a CONS in }$L_{2}\left(  [0,T]\right)  .${ Denote
by }$\mathcal{B}${ the class of all CONS in }$\mathbf{H}${ of the form
}$\left\{  e_{k}=e_{k}(s)=m_{k_{1}}(s)\ell_{k_{2}}\right\}  ${ and such that
for each }$k${, }$\sup_{0\leq s\leq T}|m_{k}(s)|<\infty.${ Obviously, for each
}$k${, }$\sup_{0\leq s\leq T}|e_{k}(s)|_{Y}<\infty.${ Let us fix a CONS
}$\mathbf{b}=\left\{  e_{k},k\geq1\right\}  \in${ }$\mathcal{B}${.}

{Let }$(\Omega,\mathcal{F}^{W},\mathbf{P})${ be a probability space with a
cylindrical Brownian motion }$W_{t}$ {in }$Y$ and $\mathcal{F}^{W}$ be the
$\sigma$-algebra generated by $W$. Let $\mathbb{F}^{W}${\ be the right
continuous filtration of }$\sigma${-algebras }$(\mathcal{F}_{t}^{W})_{t\geq0}$
generated by $W_{t}${. All the }$\sigma$-{algebras are assumed to be
}$\mathbf{P}$-{completed}$.${ Hence}%
\[
W_{t}=\sum_{k=1}^{\infty}w_{t}^{k}\ell_{k},
\]
where $\left\{  w_{t}^{k},k\geq1\right\}  ${ is a sequence of independent
standard one-dimensional Brownian motions in }$(\Omega,\mathcal{F}%
,\mathbf{P})${. We write }$W(e_{k})=\int_{0}^{T}e_{k}(t)dW_{t}$. For
$e_{k}=e_{k}(s)=m_{k_{1}}(s)\ell_{k_{2}},$
\[
W(e_{k})=\int_{0}^{T}e_{k}(t)dW_{t}=\int_{0}^{T}m_{k_{1}}(t)dw_{k_{2}}(t),
\]
and%
\begin{equation}
W_{t}=\sum_{k=1}^{\infty}(\int_{0}^{t}e_{k}(s)ds)W(e_{k}),0\leq t\leq T.
\label{f3}%
\end{equation}

{Let }$\alpha=\left\{  \alpha_{k},\text{ }k\geq1\right\}  ${ be a multiindex,
i.e. for every }$k,${ }$\alpha_{k}\in\mathbf{N}_{0}=\{0,1,2,\ldots\}.${ We
shall consider only such }$\alpha${ that }$|\alpha|=\sum_{k}\alpha_{k}<\infty
${, i.e., only a finite number of }$\alpha_{k}${ is non-zero, and we denote by
}$J${ \ the set of all such multiindices. For }$\alpha,\beta\in J${, we define
}%
\[
\alpha+\beta=(\alpha_{1}+\beta_{1},\alpha_{2}+\beta_{2},\ldots),\quad
\alpha!=\prod_{k\geq1}\alpha_{k}!.
\]
{By }$\varepsilon_{k}${ we denote the multi-index }$\alpha${ with }$\alpha
_{k}=1${ and }$\alpha_{j}=0${ for }$j\not =k${. }$\ ${Write}%

\begin{equation}
\alpha\left(  k\right)  \mathrm{=}\alpha-\varepsilon_{k} \label{eq:ek}%
\end{equation}

{For }$\alpha\in J,${ write }$H_{\alpha}:=\prod_{k=1}^{\infty}H_{\alpha_{k}%
}(W\left(  e_{k}\right)  ),${ where }$H_{n}${ is the }$n^{th}${ Hermite
polynomial defined by }$H_{n}(x)=(-1)^{N}\left(  d^{n}e^{-x^{2}/2}%
/dx^{n}\right)  e^{x^{2}/2}.$

{Let }$\xi_{\alpha}=H_{\alpha}/\sqrt{a!}${. }

\begin{theorem}
[Cameron and Martin \cite{CM}]{\label{th:CM} The set }$\Xi=\left\{
\xi_{\alpha}=\xi_{\alpha}(\mathbf{b}),\alpha\in J\right\}  ${ is an
orthonormal basis in }$L_{2}\left(  \Omega,\mathcal{F}^{W},\mathbf{P}\right)
$, where $\mathcal{F}^{W}$ is the $\sigma$-algebra generated by $W${. If }$E${
is a Hilbert space, }$\eta\in L_{2}\left(  \Omega,\mathcal{F}^{W}%
,\mathbf{P;}E\right)  ${ and }$\eta_{\alpha}=\mathbf{E}(\eta\xi_{\alpha})${,
then }$\eta=\sum_{\alpha\in{\mathcal{J}}}\eta_{\alpha}\xi_{\alpha}${ and
}$\mathbf{E}|\eta|_{E}^{2}=\sum_{\alpha\in{\mathcal{J}}}|\eta_{\alpha}%
|_{E}^{2}.$
\end{theorem}

The expansion $\eta=\sum_{\alpha\in{\mathcal{J}}}\eta_{\alpha}\xi_{\alpha}$ is
often referred to as Wiener chaos expansion.

\begin{remark}
\label{rem1}The basis $\xi_{\alpha},\alpha\in J,$ can be obtained by
differentiating stochastic exponent. Let $\mathcal{Z}$ be the set of all
real-valued sequences $z=\left(  z_{k}\right)  $ such that only finite number
of $z_{k}$ is not zero. For $\alpha\in J$, denote $\partial_{z}^{\alpha}%
=\Pi_{k}\partial^{\alpha_{k}}/\left(  \partial z_{k}\right)  ^{\alpha_{k}}$
and let%
\[
e_{z}=e_{z}(t)=\sum_{k}z_{k}e_{k}(t),0\leq t\leq T,
\]%
\begin{align}
p_{t}(z)  &  =p_{t}(e_{z})=p_{t}(z,\mathbf{b})=\exp\left\{  \int_{0}^{t}%
e_{z}(s)dW_{s}-\frac{1}{2}\int_{0}^{t}|e_{z}(s)|_{Y}^{2}ds\right\}
,\label{f30}\\
p(z)  &  =p_{T}(z),z\in\mathcal{Z},0\leq t\leq T.\nonumber
\end{align}
{It is a standard fact }(see, for example, \cite{MRWC}) {that }$H_{\alpha
}=\partial_{z}^{\alpha}p(z)|_{z=0},\xi_{\alpha}=H_{\alpha}/\sqrt{\alpha!}.${
Since }$p(z)${ is analytic, }it follows by (\ref{Eq:XiAlpha}),{ }%
\begin{equation}
p(z)=p(z,\mathbf{b})=\sum_{\alpha}\frac{H_{\alpha}}{\alpha!}z^{\alpha}%
=\sum_{\alpha}\frac{z^{\alpha}}{\sqrt{\alpha!}}\xi_{\alpha}. \label{f4}%
\end{equation}

\end{remark}

\subsection{{Generalized random variables and processes}}

Let $\mathbf{b\in}\mathcal{B},$ and $\xi_{\alpha}=\xi_{\alpha}(\mathbf{b}%
),\alpha\in J$. {Let }%
\begin{align*}
\mathcal{D}  &  =\mathcal{D}(\mathbf{b})\\
&  \mathcal{=}\left\{  v=\sum_{\alpha}v_{\alpha}\xi_{\alpha}:\text{ }%
v_{\alpha}\in\mathbf{R}\text{ and only finite number of }v_{\alpha}\text{ are
not zero}\right\}  .
\end{align*}

\begin{definition}
{A generalized }$\mathcal{D}${-random variable with values in a convex
topological vector (linear) space }$E${ }with Borel $\sigma$-algebra{ is a
formal series }$u=\sum_{\alpha}u_{\alpha}\xi_{\alpha}${, where }$u_{\alpha}\in
E,\xi_{\alpha}=\xi_{\alpha}(\mathbf{b}),${ and }$\mathbf{b}=\left\{
e_{k},k\geq1\right\}  \in\mathcal{B}${ is a CONS in }$\mathbf{H}%
=L_{2}([0,T],Y).$
\end{definition}

{Denote the vector space of all generalized }$\mathcal{D}${-random variables
by }$\mathcal{D}^{\prime}=\mathcal{D}^{\prime}(\mathbf{b})=\mathcal{D}%
^{\prime}(\mathbf{b;}E)${. The elements of }$\mathcal{D}${ are the test random
variables for }$\mathcal{D}^{\prime}.${ We define the action of a generalized
random variable }$u${ on the test random variable }$v${ by }$\left\langle
u,v\right\rangle =\sum_{\alpha}v_{\alpha}u_{\alpha}.$

{For a sequence }$u^{n}\in\mathcal{D}^{\prime}${ and }$u\in\mathcal{D}%
^{\prime}${, we say that }$u^{n}\rightarrow u${, if for every }$v\in
\mathcal{D},$ $\left\langle u,v^{n}\right\rangle \rightarrow\left\langle
u,v\right\rangle .$ {This implies that }$u^{n}=\sum_{\alpha}u_{\alpha}^{n}%
\xi_{\alpha}\rightarrow u=\sum_{\alpha}u_{\alpha}\xi_{\alpha}${ if and only if
}$u_{\alpha}^{n}\rightarrow u_{\alpha}${ as }$n\rightarrow\infty${ for all
}$\alpha${.}

\begin{remark}
\label{rem11}Obviously, if $u=\sum_{\alpha}u_{\alpha}\xi_{\alpha}%
\in\mathcal{D}^{\prime}(\mathbf{b};E),~F$ is a vector space and
$f:E\rightarrow F$ is a linear map, then%
\[
f(u)=\sum_{\alpha}f(u_{\alpha})\xi_{\alpha}\in\mathcal{D}^{\prime}%
(\mathbf{b};F).
\]

\end{remark}

\begin{definition}
{An\ }$E${-valued generalized }$\mathcal{D}${- process\ }$u(t)${ in }$[0,T]${
is a }$\mathcal{D}^{\prime}(\mathbf{b};E)${-valued function on }$[0,T]${ such
that for each }$t\in\lbrack0,T]$%
\[
u(t)=\sum_{\alpha}u_{\alpha}(t)\xi_{\alpha}\in\mathcal{D}^{\prime}%
(\mathbf{b};E);
\]
and $u_{\alpha}(t)$ are deterministic measurable $E$-valued functions on
$[0,T]$. We denote the linear space of all such processes by{ }$\mathcal{D}%
^{\prime}(\mathbf{b;}[0,T],E).$ If $E$ is a topological vector space and a
generalized $\mathcal{D}${-process} $u(t)$ is continuous we write $u\in$
$C\mathcal{D}^{\prime}([0,T],\mathbf{b},E)$ (note that $u(t)$ is continuous if
and only if all coefficient functions $u_{\alpha}$ are continuous in $E$.

{If there is no room for confusion, we will often say }$\mathcal{D}${-process
(}$\mathcal{D}${-random variable) instead of generalized }$\mathcal{D}%
${-process (generalized }$\mathcal{D}${-random variable).}
\end{definition}

{If }$E${ is a normed vector space, we denote}%
\begin{align*}
&  L_{1}(\mathcal{D}^{\prime}\mathcal{(}\mathbf{b};[0,T],E))\\
&  =\{u(t)=\sum_{\alpha}u_{\alpha}(t)\xi_{\alpha}\in\mathcal{D}^{\prime
}\mathcal{(}\mathbf{b};[0,T],E):\int_{0}^{T}|u_{\alpha}(t)|_{E}dt<\infty
\text{,~}\alpha\in J\}.
\end{align*}
{For }$u(t)=\sum_{\alpha}u_{\alpha}(t)\xi_{\alpha}\in L_{1}(\mathcal{D}%
^{\prime}(\mathbf{b};[0,T],E))${ we define }$\int_{0}^{t}u(s)ds,0\leq t\leq
T,${ in }$\mathcal{D}^{\prime}(\mathbf{b};[0,T],E)${ by}%
\[
\int_{0}^{t}u(s)ds=\sum_{\alpha}\left(  \int_{0}^{t}u_{\alpha}(s)ds\right)
\xi_{\alpha},0\leq t\leq T.
\]

{If }$u(t)=\sum_{\alpha}u_{\alpha}(t)\xi_{\alpha}\in\mathcal{D}^{\prime
}\mathcal{(}\mathbf{b};[0,T],E),$ then $u(t)$ is differentiable in $t$ if and
only if $u_{\alpha}(t)${ are differentiable in }$t\,$. In that case{, }%
\[
\frac{d}{dt}u(t)=\dot{u}(t)=\sum_{\alpha}\dot{u}_{\alpha}(t)\xi_{\alpha}%
\in\mathcal{D}^{\prime}\mathcal{(}[0,T],\mathbf{b},E).
\]

\begin{example}
{\label{ex:CWP}A cylindrical Wiener process }$W_{t},0\leq t\leq T,${ in a
Hilbert space }$Y,${ and its derivative }$dW_{t}/dt=\dot{W}_{t}${ are
generalized }$Y${-valued stochastic processes. Indeed, by (\ref{f3}),}%
\[
W_{t}=\sum_{k}\int_{0}^{t}e_{k}(s)ds\xi_{\varepsilon_{k}},0\leq t\leq T,
\]
{and }$W_{t}=\int_{0}^{t}\dot{W}_{s}ds,$ {where }$\dot{W}_{t}=\sum_{k}%
e_{k}(t)\xi_{\varepsilon_{k}},0\leq t\leq T.$
\end{example}

\subsubsection{{Wick Product and Skorokhod Integral\label{WickSkor}}}

\begin{definition}
{\label{def:WP} For }$\xi_{\alpha}${, }$\xi_{\beta}${ from }$\Xi${, define the
Wick product }%
\begin{equation}
\xi_{\alpha}\lozenge\xi_{\beta}:=\sqrt{\left(  \frac{(\alpha+\beta)!}%
{\alpha!\beta!}\right)  }\xi_{\alpha+\beta}. \label{eq:WP}%
\end{equation}

\end{definition}

{In particular, taking in (\ref{eq:WP}) }$\alpha=k\varepsilon_{i}${ and
}$\beta=n\varepsilon_{i}${ we get }%
\begin{equation}
H_{k}(\xi_{i})\lozenge H_{n}(\xi_{i})=H_{k+n}(\xi_{i}). \label{eq:WP-hp}%
\end{equation}
{For a Hilbert space }$E${ and arbitrary }$v=\sum_{\alpha}v_{\alpha}%
\xi_{\alpha}${ and }$u=\sum_{\alpha}u_{\alpha}\xi_{\alpha}${ in }%
$\mathcal{D}^{\prime}(\mathbf{b};E)${, we define their Wick product as a
}$\mathcal{D}${-generalized real valued random variable given by}%
\begin{equation}
v\lozenge u=\sum_{\alpha}\sum_{\beta\leq\alpha}(u_{\beta},v_{\alpha-\beta
})_{E}\sqrt{\frac{\alpha!}{\beta!(\alpha-\beta)!}}\xi_{\alpha}\in
\mathcal{D}^{\prime}\mathcal{(}\mathbf{b};\mathbf{R}). \label{eq:WPDprime}%
\end{equation}

\begin{definition}
\label{de1}{Skorokhod integral (Maliavin divergence operator) of\ }$v\in$%
{\ }$L_{1}(\mathcal{D}^{\prime}\mathcal{(}\mathbf{b};[0,T],Y))${ is a
generalized random variable }$($element of $\mathcal{D}^{\prime}\left(
\mathbf{b};\mathbf{R}\right)  $){ such that }%
\[
\delta(v)=\int_{0}^{T}v(s)dW_{s}=\sum_{\alpha}\delta(v)_{\alpha}\xi_{\alpha},
\]
{with}%
\begin{equation}
\delta(v)_{\alpha}=\sum_{k}\sqrt{\alpha_{k}}\int_{0}^{T}\left(  v_{\alpha
(k)}(t),e_{k}(t)\right)  _{Y}dt. \label{f1}%
\end{equation}
{and }$\alpha(k)${ is given by (\ref{eq:ek}). }
\end{definition}

If $v\in${\ }$L_{1}(\mathcal{D}^{\prime}\mathcal{(}\mathbf{b};[0,T],Y))$, then
$\delta_{t}(v)=\int_{0}^{t}v(s)dW_{s}=\delta\left(  v1_{[0,t]}\right)  ,0\leq
t\leq T,${ is a process in }$\mathcal{D}^{\prime}\left(  \mathbf{b}%
;[0,T],Y\right)  $. We have%
\[
\delta_{t}(v)_{\alpha}=\sum_{k}\sqrt{\alpha_{k}}\int_{0}^{t}\left(
v_{\alpha(k)}(s),e_{k}(s)\right)  _{Y}ds.
\]
{ Since }$\dot{W}_{t}=\sum_{k}e_{k}(t)\xi_{\varepsilon_{k}},${ it follows by
(\ref{eq:WPDprime}) that}%
\[
v_{t}\lozenge\dot{W}_{t}=\sum_{\alpha}\sum_{k}(v_{\alpha(k)}(t),e_{k}%
(t))_{Y}\sqrt{\alpha_{k}}\xi_{\alpha},
\]
and%
\[
\delta\left(  v\right)  =\int_{0}^{T}v_{t}\lozenge\dot{W}_{t}~dt,\;\delta
_{t}(v)=\int_{0}^{t}v(s)\lozenge\dot{W}_{s}ds,0\leq t\leq T.
\]

\begin{remark}
\label{rem2}Skorokhod integral is an extension of the It\^{o}
integral\footnote{Of course, this statement is well known. However, the proof
given here is short and straightforward.}; (\ref{form1}) below motivates the
definition of the Skorokhod integral.
\end{remark}

If $u(t)=\sum_{\alpha}u_{\alpha}(t)\xi_{\alpha}$ is $\mathbb{F}^{W}$-adapted
$Y$-valued such that%
\[
\mathbf{E}\int_{0}^{T}|u(t)|_{H}^{2}dt<\infty,
\]
then $v=\int_{0}^{T}u(t)dW(t)=\sum_{\alpha}v_{\alpha}\xi_{\alpha}$ is square
integrable. By Ito formula for the product of $\int_{0}^{t}u(s)dW(s)$ and
stochastic exponent $p_{t}(z)$ from Remark \ref{rem1}, we obtain%
\begin{align*}
\mathbf{E}vp(z)  &  =\mathbf{E}\int_{0}^{T}u(t)dW(t)p_{T}(z)=\int_{0}%
^{T}\mathbf{E[}p_{t}(z)(u(t),e_{z}(t))_{Y}]dt\\
&  =\int_{0}^{T}\mathbf{E[}p(z)(u(t),e_{z}(t))_{Y}]dt,z\in\mathcal{Z}.
\end{align*}
So,%
\begin{align}
\frac{\partial^{|\alpha|}\mathbf{E}vp(z)}{\partial z^{\alpha}}  &  =\sum
_{k}\alpha_{k}\int_{0}^{T}\frac{\partial^{|\alpha(k)|}}{\partial z^{\alpha
(k)}}(\mathbf{E}p(z)u(t),e_{k}(t))_{Y}]dt\mathrm{,}\nonumber\\
v_{\alpha}  &  =(\sqrt{\alpha!})^{-1}\frac{\partial^{|\alpha|}\mathbf{E}%
vp(z)}{\partial z^{\alpha}}|_{z=0}=\sum_{k}\int_{0}^{T}\sqrt{\alpha_{k}%
}(u_{\alpha(k)}(t),e_{k}(t))_{Y}]dt. \label{form1}%
\end{align}
Comparing (\ref{form1}) and (\ref{f1}), we see that Ito and Skorokhod
integrals are equal in this case.

\section{{ \ \ }Wick product{ Navier-Stokes Equation}}

{For }$T>r\geq0,${ let us consider the following Navier-Stokes equation:}%
\begin{align}
\partial_{t}\mathbf{u}\left(  t,x\right)   &  =\partial_{i}\left(
a^{\,ij}\left(  t,x\right)  \partial_{j}\mathbf{u}\left(  t,x\right)  \right)
+b^{i}(t,x)\partial_{i}\mathbf{u}(t,x)\label{nse}\\
&  -u^{k}\left(  t,x\right)  \lozenge\partial_{k}\mathbf{u}%
(t,x)+\mathbf{\nabla}P\left(  t,x\right)  +\mathbf{f}\left(  t,x\right)
\nonumber\\
&  \lbrack\sigma^{i}(t,x)\partial_{i}\mathbf{u}\left(  t,x\right)
+\mathbf{g}\left(  t,x\right)  -\mathbf{\nabla}\tilde{P}\left(  t,x\right)
]\lozenge\dot{W}_{t},\nonumber\\
\text{ }\mathbf{u}\left(  r,x\right)   &  =\boldsymbol{w}\left(  x\right)
,\operatorname{div}\text{\thinspace}\mathbf{u}=0.\nonumber
\end{align}
{The unknowns in the equation (\ref{nse}) are the functions }$\mathbf{u}%
=\left(  u^{l}\right)  _{1\leq l\leq d},\,P,\,\tilde{P}${. It is assumed that
}$a^{ij},b^{i},\mathbf{f}=\left(  f^{i}\right)  ,${ are measurable
deterministic functions on }$[0,\infty)\times\mathbf{R}^{d}${, and the matrix
}$\left(  a^{ij}\right)  ${ is symmetric.\ Let us assume also that }%
$\sigma^{i},\mathbf{g}=\left(  g^{i}\right)  ${ be }$Y${-valued measurable
deterministic functions on }$[0,\infty)\times\mathbf{R}^{d}.$ Let $\mathbf{w}$
be a random initial velocity field.

{In addition, we will need the following assumptions.}

\textbf{A1}{. For all }$t\geq0,x\in\mathbf{R}^{d}{,}\lambda\in\mathbf{R}^{d}%
,${ }%
\[
K|\lambda|^{2}\geq a^{ij}(t,x)\lambda_{i}\lambda_{j}\geq\delta|\lambda|^{2},
\]
{where }$K,\delta${ are fixed strictly positive constants.}

\textbf{A2}{. For all }$t\geq0,x,${ }%
\[
\max_{|\alpha|\leq2}|\partial^{\alpha}a^{ij}(t,x)|+\max_{|\alpha|\leq
1}(|\partial^{\alpha}\,b^{i}(t,x)|+|\partial^{\alpha}\,\sigma^{i}%
(t,x)|_{Y})\leq K.
\]

\textbf{A3}{. The functions}\textbf{ }$\mathbf{f}(t,x)${ and }$\mathbf{g}%
(t,x)${ are measurable deterministic, }$p>d${, and for all }$t>0,$%
\[
\int_{0}^{t}\sum_{l=2,p}[|\mathbf{f(}r\mathbf{)|}_{1,l}^{l}+||\mathbf{g(}%
r\mathbf{)||}_{1,l}^{2l}]dr<\infty
\]
(recall $|\mathbf{f(}r\mathbf{)|}_{1,l},||\mathbf{g(}r\mathbf{)||}_{1,l}$ are
$\mathbb{H}_{p}^{1}(\mathbf{R}^{d})$ and $\mathbb{H}_{p}^{1}(\mathbf{R}%
^{d},Y)$-norms respectively).

{We will seek a solution to (\ref{nse}) in the form }%
\[
\mathbf{u}(t)\mathbf{=}\sum_{\alpha}\mathbf{u}_{\alpha}(t)\xi_{\alpha}%
\in\mathcal{D}^{\prime}(\mathbf{b};[0,T]\mathbf{,}\mathbb{H}_{p}^{2}),p\geq2.
\]
{In this case, denoting by }$\mathcal{P}(\mathbf{v})${ the solenoidal
projection of the vector field }$\mathbf{v}${, we can rewrite (\ref{nse}) in
the following equivalent form:\ }%
\begin{align}
\partial_{t}\mathbf{u}\left(  t\right)   &  =\mathcal{P}[\partial_{i}\left(
a^{\,ij}\left(  t\right)  \partial_{j}\mathbf{u}\left(  t\right)  \right)
+b^{i}(t)\partial_{i}\mathbf{u}(t)\label{nse0}\\
&  -u^{k}\left(  t\right)  \partial_{k}\mathbf{u}(t)+\mathbf{f}\left(
t\right)  ]+\mathcal{P}[\sigma^{i}(t)\partial_{i}\mathbf{u}\left(  t\right)
+\mathbf{g}\left(  t\right)  ]\lozenge\dot{W}_{t},\nonumber\\
\mathbf{u}\left(  r\right)   &  =\boldsymbol{w,}\operatorname{div}%
\text{\thinspace}\mathbf{u}(t)=0,t\in\lbrack r,T].\nonumber
\end{align}
If $\mathbf{\eta}=\sum_{\alpha}\mathbf{\eta}_{\alpha}\xi_{\alpha}$ with
$\mathbf{\eta}_{\alpha}\in\mathbb{H}_{p}^{k}$, then (see Remark \ref{rem11})
$\mathcal{P}(\mathbf{\eta})=$ $\sum_{\alpha}\mathcal{P}(\mathbf{\eta}_{\alpha
})\xi_{\alpha}$.

We start our analysis of equation (\ref{nse0}) by introducing the definition
of a solution in the "weak sense".

\begin{definition}
{\label{defs}We say that a generalized }$\mathcal{D}${-process }%
$\mathbf{u}(t)=\sum_{\alpha}\mathbf{u}_{\alpha}(t)\xi_{\alpha}\in
C\mathcal{D}^{\prime}(\mathbf{b};[r,T],\mathbb{H}_{p}^{k})${ is }$\mathcal{D}%
${-}$\mathbb{H}_{p}^{k}${ solution of equation (\ref{nse}) in }$[r,T]${, if
the equality }%
\begin{align}
\mathbf{u}(t)  &  =\mathbf{w}+\int_{r}^{t}\mathcal{P}[-u^{i}\left(  s\right)
\lozenge\partial_{i}\mathbf{u}\left(  s\right)  +\partial_{i}(a^{ij}%
(s)\partial_{j}\mathbf{u}\left(  s\right)  \mathbf{)}\label{nse1}\\
&  \mathbf{+}b^{i}(s)\mathbf{\partial}_{i}\mathbf{\mathbf{u}}(s)+\mathbf{f}%
(s)]ds\nonumber\\
&  \int_{r}^{t}\mathcal{P}[\sigma^{k}(s)\partial_{k}\mathbf{u}\left(
s\right)  +\mathbf{g}(s)]\lozenge\dot{W}_{s}ds\nonumber
\end{align}
{holds in }$\mathcal{D}(\mathbf{b};\mathbb{H}_{p}^{k-2}(\mathbf{R}^{d}))${ for
every }$r\leq t\leq T${. If an }$\mathcal{D}$-$\mathbb{H}_{p}^{k}${-solution
in }$[r,T]${ is also }$\mathcal{D}$-$\mathbb{H}_{q}^{k^{\prime}}${-solution in
}$[r,T]${, we call it }$\mathcal{D}$-$\mathbb{H}_{p}^{k}\cap\mathbb{H}%
_{q}^{k^{\prime}}${-solution in }$[r,T]${. }

In the future, we simply say $\mathcal{D}$-solution if there is no risk of confusion.
\end{definition}

\begin{remark}
\label{rem12}1. {Assume A1-A3 hold, }$p\geq2,\mathbf{w=}\sum_{\alpha
}\mathbf{w}_{\alpha}\xi_{\alpha}\in\mathcal{D}^{\prime}(\mathbf{b}%
,\mathbb{H}_{p}^{k})${. Applying Remark \ref{rem11} and definition of the Wick
product we see that }$\mathbf{u}(t)=\sum_{\alpha}\mathbf{u}_{\alpha}%
(t)\xi_{\alpha}\in C\mathcal{D}^{\prime}([r,T],\mathbf{b},\mathbb{H}_{p}^{k}%
)${ is an }$\mathcal{D}${-}$\mathbb{H}_{p}^{k}${ solution in }$[r,T]${ if and
only if for each }$\alpha$, $\mathbf{u}_{\alpha}\in C([0,T],${ $\mathbb{H}%
_{p}^{k})$} and for $t\in\lbrack r,T]$ the following equality holds in
$\mathbb{H}_{p}^{k-2}:$%
\begin{align}
\mathbf{u}_{\alpha}\left(  t\right)   &  =\mathbf{w}_{\alpha}+\int_{r}%
^{t}\mathcal{P}\{\partial_{i}\left(  a^{\,ij}\left(  s\right)  \partial
_{j}\mathbf{u}_{\alpha}\left(  s\right)  \right)  +b^{i}(s)\partial
_{i}\mathbf{u}_{\alpha}(s)\label{pr1}\\
&  -\sum_{\gamma\leq\alpha}\sqrt{\binom{\alpha}{\gamma}}u_{\alpha-\gamma}%
^{k}\left(  s\right)  \partial_{k}\mathbf{u}_{\gamma}(s)+\mathbf{f}\left(
s\right)  1_{\alpha=0}\nonumber\\
&  +\sum_{k}\sqrt{\alpha_{k}}[\left(  \sigma^{i}(s),e_{k}(s)\right)
_{Y}\partial_{i}\mathbf{u}_{\alpha(k)}\left(  s\right) \nonumber\\
&  +\left(  \mathbf{g}\left(  s\right)  ,e_{k}(s)\right)  _{Y}1_{|\alpha
|=1}]ds.\nonumber
\end{align}

2. {If }$\alpha=0${, the zero term }$\mathbf{u}_{\alpha}(t,x)=\mathbf{u}%
_{0}(t,x)${ of an }$\mathcal{D}${-}$\mathbb{H}_{p}^{k}${ solution in }$[r,T]${
satisfies Navier-Stokes equation: }%
\begin{align}
\mathbf{u}_{0}\left(  t\right)   &  =\mathbf{w}_{0}+\int_{r}^{t}%
\mathcal{P}[\partial_{i}\left(  a^{\,ij}\left(  s\right)  \partial
_{j}\mathbf{u}_{0}\left(  s\right)  \right)  +b^{i}(s)\partial_{i}%
\mathbf{u}_{0}(s)\label{pf1}\\
&  -u_{0}^{k}\left(  s\right)  \partial_{k}\mathbf{u}_{0}(s)+\mathbf{f}\left(
s\right)  ]ds.\nonumber
\end{align}

{For the remaining components we have to solve Stokes equations. For }%
$|\alpha|\geq1,$ we can rewrite (\ref{pr1}) as%
\begin{align}
\mathbf{u}_{\alpha}\left(  t\right)   &  =\mathbf{w}_{\alpha}+\int_{r}%
^{t}\mathcal{P}[\partial_{i}\left(  a^{\,ij}\left(  s\right)  \partial
_{j}\mathbf{u}_{\alpha}\left(  s\right)  \right)  +\mathbf{F}_{\alpha
}(s)\label{pf4}\\
&  +[b^{i}(s)-u_{0}^{i}(s)]\partial_{i}\mathbf{u}_{\alpha}(s)-u_{\alpha}%
^{k}\left(  s\right)  \partial_{k}\mathbf{u}_{0}(s)]ds,\nonumber
\end{align}
with
\begin{align}
\mathbf{F}_{\alpha}(s)  &  =\sum_{\gamma\leq\alpha,|\alpha|-1\geq|\gamma
|\geq1}\sqrt{\binom{\alpha}{\gamma}}u_{\alpha-\gamma}^{k}\left(  s\right)
\partial_{k}\mathbf{u}_{\gamma}(s)\label{pf5}\\
&  +\sum_{k}\sqrt{\alpha_{k}}[\left(  \sigma^{i}(s),e_{k}(s)\right)
_{Y}\partial_{i}\mathbf{u}_{\alpha(k)}\left(  s\right)  +\left(
\mathbf{g}\left(  s\right)  ,e_{k}(s)\right)  _{Y}1_{|\alpha|=1}].\nonumber
\end{align}

{Since for }$|\alpha|\geq1,${ }$\mathbf{E}\xi_{\alpha}=0${, the equation
(\ref{nse0}) (or (\ref{nse1})) can be regarded as a random perturbation of the
deterministic Navier-Stokes equation (\ref{pf1}).}
\end{remark}

\begin{lemma}
\label{leme}{Let A1-A3 hold, }$|\mathbf{w}_{0}|_{2,p}+|\mathbf{w}_{0}%
|_{2,2}<\infty.$ {Then there is }$T_{1}>0${ and a unique }$\mathbf{u}_{0}\in
C([0,T_{1}),\mathbb{H}_{2}^{2}\cap\mathbb{H}_{p}^{2})$ solving {(\ref{pf1}) in
}$[0,T_{1}).$
\end{lemma}

\begin{proof}
{According to Theorem 3 in \cite{mr22}, there is }$T_{1}>0${ and a unique
}$\mathbf{u}_{0}\in C\left(  [0,T_{1}\right)  ,\mathbb{H}_{p}^{1}%
\cap\mathbb{H}_{2}^{1})$ {such that for each }$t<T_{1}${ }%
\[
\sup_{0\leq s\leq t}|\mathbf{u}_{0}(s)|_{1,l}^{l}+\int_{r}^{t}|\partial
^{2}\mathbf{u}_{0}(r)|_{l}^{l}dr<\infty,l=2,p,
\]
and {(\ref{pf1}) holds in }$\mathbb{H}_{l}^{-1},l=2,p.${ By Sobolev embedding
theorem, for all }$t<T_{1},$%
\[
\sup_{x,s\leq t}|\mathbf{u}_{0}(s,x)|+\int_{0}^{t}\sup_{x}|\nabla
\mathbf{u}_{0}(s,x)|^{p}ds<\infty.
\]
{and there is a constant }$C$ such that {for all }$s\in\lbrack0,T_{1}),$%
\[
|u_{0}^{k}(s)\partial_{k}\mathbf{u}_{0}(s)|_{1,p}\leq C|u_{0}^{k}%
(s)|_{1,p}|\partial_{k}\mathbf{u}_{0}(s)|_{1,p}.
\]
So, for each $t<T_{1},$%
\begin{align*}
&  \int_{0}^{t}|u_{0}^{k}(s)\partial_{k}\mathbf{u}_{0}(s)|_{1,p}^{p}ds\leq
C\int_{0}^{t}|u_{0}^{k}(s)|_{1,p}^{p}|\partial_{k}\mathbf{u}_{0}(s)|_{1,p}%
^{p}ds\\
&  \leq C\sup_{s\leq t}|\mathbf{u}_{0}(s)|_{1,p}\int_{0}^{t}|\partial
_{k}\mathbf{u}_{0}(s)|_{1,p}^{p}ds<\infty.
\end{align*}
Also,
\[
\int_{0}^{t}|u_{0}^{k}(s)\partial_{k}\mathbf{u}_{0}(s)|_{1,2}^{2}ds<\infty.
\]
Indeed,
\[
\int_{0}^{t}|u_{0}^{k}(s)\partial_{k}\mathbf{u}_{0}(s)|_{2}^{2}ds\leq
\sup_{s\leq t,x}|u_{0}^{k}(s,x)|^{2}\int_{0}^{t}|\nabla\mathbf{u}_{0}%
(s)|_{2}^{2}ds<\infty,
\]
and also
\begin{align*}
&  \int_{0}^{t}|\nabla\left(  u_{0}^{k}(s)\partial_{k}\mathbf{u}%
_{0}(s)\right)  |_{2}^{2}ds=I_{1}+I_{2}\\
&  =\int_{0}^{t}|\nabla\mathbf{u}_{0}(s)|_{4}^{4}ds+\int_{0}^{t}|u_{0}%
^{k}(s)\partial_{k}\nabla\mathbf{u}_{0}(s)|_{2}^{2}ds]
\end{align*}
with%
\begin{align*}
I_{1}  &  \leq\int_{0}^{t}\sup_{x}|\nabla\mathbf{u}_{0}(s,x)|^{2}%
|\nabla\mathbf{u}_{0}(s)|_{2}^{2}ds\\
&  \leq\sup_{s\leq t}|\nabla\mathbf{u}_{0}(s)|_{2}^{2}\int_{0}^{t}%
|\mathbf{u}_{0}(s)|_{2,p}^{2}<\infty,
\end{align*}
and%
\[
I_{2}\leq\sup_{s\leq t,x}|\mathbf{u}_{0}(s,x)|^{2}\int_{0}^{t}|\partial
^{2}\mathbf{u}_{0}(s)|_{2}^{2}ds<\infty.
\]

By {Proposition \ref{aprop1} in Appendix II, }${\mathbf{u}}_{0}\in C\left(
[0,T],\mathbb{H}_{p}^{2}\cap\mathbb{H}_{2}^{2}\right)  $ for every $T<T_{1}$
and (\ref{pf1}) holds in $\mathbb{L}_{l},l=2,p$.
\end{proof}

Now we fix an arbitrary $T<T_{1}$ ($T_{1}$ comes from Lemma \ref{leme}) and
{prove the existence and uniqueness of }$\mathcal{D}${-solutions to
(\ref{nse0}) in }$[r,T],r<T${. }

\begin{lemma}
{\label{l1}Assume that A1-A3 hold, }$\mathbf{b}\in\mathcal{B},\mathbf{w}%
=\sum_{\alpha}\mathbf{w}_{\alpha}\xi_{\alpha}\in\mathcal{D}^{\prime
}\mathbf{(b};\mathbb{H}_{p}^{2}\cap\mathbb{H}_{2}^{2}).$

{Then for each }$r\leq T<T_{1}${ there is a unique }$\mathcal{D}$%
{-}$\mathbb{H}_{p}^{2}\cap\mathbb{H}_{2}^{2}${-solution }$\mathbf{u}%
(t)=\sum_{\alpha}\mathbf{u}_{\alpha}(t)\xi_{\alpha}\in C\mathcal{D}^{\prime
}\left(  \mathbf{b;}[0,T],\mathbb{H}_{p}^{2}\cap\mathbb{H}_{2}^{2}\right)  $
{of (\ref{nse}) in }$[r,T].$

Equivalently, for each $\alpha,$ $\mathbf{u}_{\alpha}\in C\left(
[0,T],\mathbb{H}_{p}^{2}\cap\mathbb{H}_{2}^{2}\right)  $ and (\ref{pf1}%
)-(\ref{pf4}) hold in $\mathbb{L}_{l},l=2,p.$
\end{lemma}

\begin{proof}
According to Remark \ref{rem12}, it suffice to prove the existence and
uniqueness of a solution for the deterministic system (\ref{pr1}). {For
}$\alpha=0,${ the existence and uniqueness of a solution to (\ref{pf1})
follows from Lemma \ref{leme}. }We proceed by induction. {Assume there are
unique }$\mathbf{u}_{\alpha}\in C([0,T],\mathbb{H}_{p}^{2}\cap\mathbb{H}%
_{2}^{2}),|\alpha|\leq n${,} such that (\ref{pr1}) holds in $\mathbb{L}%
_{l},l=2,p.$ {By Sobolev embedding theorem, it implies that }%
\begin{equation}
\sup_{x,r\leq s\leq T}|\mathbf{u}_{\alpha}(s,x)|+\int_{r}^{T}\sup_{x}%
|\partial\mathbf{u}_{\alpha}(s,x)|^{p}ds<\infty, \label{es3'}%
\end{equation}
{if }$|\alpha|\leq n${. Then for }$|\alpha|=n+1,$ the equation (\ref{pf4}) for
$\mathbf{u}_{\alpha}$ is Stokes and it is readily checked (see (\ref{pf4}))
that
\[
\int_{r}^{T}|\mathbf{F}_{\alpha}(s)|_{1,l}^{l}ds<\infty,l=2,p.
\]
According to Proposition \ref{aprop1}, there is a unique $\mathbf{u}_{\alpha
}\in C\left(  [0,T],\mathbb{H}_{p}^{2}\cap\mathbb{H}_{2}^{2}\right)  $ so that
(\ref{pf4}) holds in $\mathbb{L}_{l},l=2,p.$
\end{proof}

{Because of the uniqueness, the }$\mathcal{D}${-solution has a restarting
property. More specifically, the following statement holds.}

\begin{corollary}
{\label{cornew1} Assume that A1-A3 hold, }$\mathbf{b}\in\mathcal{B}%
,\mathbf{w}=\sum_{\alpha}\mathbf{w}_{\alpha}\xi_{\alpha}\in\mathcal{D}%
^{\prime}\mathbf{(b};\mathbb{H}_{p}^{2}\cap\mathbb{H}_{2}^{2}).$ {Let
}$\mathbf{u}^{r,\mathbf{w}}(t)${ be the unique }$\mathcal{D}${-}%
$\mathbb{H}_{p}^{2}\cap\mathbb{H}_{2}^{2}$ {solution to (\ref{nse})in
}$[r,T],T<T_{1},${ starting at }$\mathbf{w}$. Let $r\leq r^{\prime}\leq t\leq
T.${ Then}%
\[
\mathbf{u}^{r,\mathbf{w}}(t)=\mathbf{u}^{r^{\prime},\mathbf{u(}r^{\prime}%
)}(t).
\]

\end{corollary}

\begin{proof}
{Indeed for }$u(t)=u^{r,\mathbf{w}}(t),${ and }$r\leq r^{\prime}\leq t\leq
T,${ we have for }$t\in\lbrack r,T]$%
\begin{align*}
\mathbf{u}\left(  t\right)   &  =\mathbf{u}(r^{\prime})+\int_{r^{\prime}}%
^{t}\mathcal{P}\{\partial_{i}\left(  a^{\,ij}\left(  s\right)  \partial
_{j}\mathbf{u}\left(  s\right)  \right)  +b^{i}(s)\partial_{i}\mathbf{u}(s)\\
&  -u^{k}\left(  s\right)  \lozenge\partial_{k}\mathbf{u}(s)+\mathbf{f}\left(
s\right) \\
&  +[\sigma^{i}(s)\partial_{i}\mathbf{u}\left(  s\right)  +\mathbf{g}\left(
s\right)  -\mathbf{\nabla}\tilde{P}\left(  s\right)  ]\lozenge\dot{W}_{s}\}ds,
\end{align*}
{and the statement follows by Lemma \ref{l1}.}
\end{proof}

\subsection{{\bigskip Rescaling and approximation of the generalized
solution}}

{To begin with, we will derive more precise estimates for }$\mathcal{D}$%
{-}$\mathbb{H}_{p}^{2}\cap\mathbb{H}_{2}^{2}${ solutions of equation
(\ref{nse0}). One could hardly expect that the }$\mathcal{D}${-}%
$\mathbb{H}_{p}^{2}\cap\mathbb{H}_{2}^{2}${ solution of unbiased Navier-Stokes
equation has finite variance, i.e. }$\sum_{\alpha}\left\vert \mathbf{u}%
_{\alpha}(t)\right\vert ^{2}<\infty.${ However, in this subsection we will
show that the solution could be \ obtained as the limit of square integrable
solutions of the equations }rescaled in a special way using a second
quantization operator{ (see \cite{Simon}, and Appendix I, 5.1).}

Fix $\mathbf{b=}\left\{  e_{n}\right\}  \mathbf{\in}\mathcal{B}$ and \ define
an unbounded operator
\[
Ae_{k}=2ke_{k},k\geq1.
\]
Obviously, the projective limit of the domains $\mathbf{H}_{n}\subseteq
\mathbf{H}=L_{2}\left(  [0,T],Y\right)  $ of $A^{n}$ with the norm%
\[
||y||_{\mathbf{H}_{n}}=||A^{n}y||=\left(  \sum_{k}\left(  2k\right)
^{2n}y_{k}^{2}\right)  ^{1/2},y=\sum_{k}y_{k}e_{k}\in\mathbf{H}_{n},
\]
is a nuclear space denoted $\mathcal{N}=\mathcal{N}(\mathbf{b})$. For
$n\in\mathbf{N}$, let $\mathbf{H}_{-n}$ be the completion of $\mathbf{H}$ with
respect to the norm%
\[
||y||_{\mathbf{H}_{-n}}=||A^{-n}y||=(\sum_{k}\left(  2k\right)  ^{-2n}%
y_{k}^{2})^{1/2},y=\sum_{k}y_{k}e_{k}\in\mathbf{H}_{-n}.
\]
The inductive limit $\mathcal{N}^{\prime}=\mathcal{N}^{\prime}(\mathbf{b)}$ of
$\mathbf{H}_{-n}$ is the dual of $\mathcal{N}$. For a Banach space $E$ and
$\rho\in\lbrack0,1]$, let $\mathcal{S}_{\rho}(E)=\mathcal{S}_{\rho}%
(\mathbf{b};E)$ be the space of all $\eta=\sum_{\alpha}a_{\alpha}\xi_{\alpha
}\in\mathcal{D}^{\prime}(\mathbf{b},E\mathbf{)}$ such that
\[
||\eta||_{\mathcal{S}_{\rho,q}}=||\eta||_{\rho,q}=(\sum_{\alpha}%
(|\alpha|!)^{\rho}(2\mathbf{N})^{2q\alpha}|a_{\alpha}|_{E}^{2})^{1/2}%
<\infty\text{ for every }q\geq0,
\]
where
\[
(2\mathbf{N})^{2q\alpha}=\prod_{k=1}^{\infty}(2k)^{2q\alpha_{k}}%
=2^{2q|\alpha|}\prod_{k=1}^{\infty}k^{2q\alpha_{k}}%
\]
(see \cite{holoks}). We consider $\mathcal{S}_{\rho}(E)=\mathcal{S}_{\rho
}(\mathbf{b},E)$ with a family of seminorms $||\eta||_{\mathcal{S}_{\rho};q}.$
Similarly, for $\rho\in\lbrack0,1]$, let $\mathcal{S}_{-\rho}(E)$ be the space
of all $\eta\in\mathcal{D}^{\prime}(\mathbf{b};E)$ such that $||\eta
||_{\mathcal{S}_{-\rho,-q}}=||\eta||_{-\rho,-q}<\infty$ for some $q>0.$ It is
dual of $\mathcal{S}_{-\rho}(E)$ if $E$ is Hilbert. If $E=\mathbf{R}$,
Kondratiev test function space $\left(  \mathcal{N}\right)  ^{1}%
=\mathcal{S}_{1}(\mathbf{R})$ and Kondratiev distribution space $\left(
\mathcal{N}\right)  ^{-1}=\mathcal{S}_{-1}(\mathbf{R})$ (see \cite{holoks},
p.39, \cite{kon3}).

We will show that the solution found in Lemma \ref{l1} belongs to the
Kondratiev space $S_{-1}(\mathbb{H}_{p}^{2}\cap\mathbb{H}_{2}^{2})=\left(
\mathcal{N}\right)  ^{-1}\otimes(\mathbb{H}_{p}^{2}\cap\mathbb{H}_{2}^{2})$.

\begin{proposition}
{\label{cle1}Let A1-A3 hold,} $\sup_{s,k}|e_{k}(s)|_{Y}<\infty$ and
$\mathbf{w}\in\mathbb{H}_{p}^{2}\cap\mathbb{H}_{2}^{2}$ be deterministic$.${
Assume that}%
\[
\mathbf{u}(t)=\sum_{\alpha}\mathbf{u}_{\alpha}(t)\xi_{\alpha}\in
C\mathcal{D}^{\prime}(\mathbf{b};[0,T],\mathbb{H}_{p}^{2}\cap\mathbb{H}%
_{2}^{2})
\]
{ } solve {(\ref{nse0}) in }$\mathbb{L}_{l},l=2,p${. Denote}%
\[
\tilde{L}_{\alpha}=\sup_{t\leq T}|\mathbf{u}_{\alpha}(t)|_{2,p}+\sup_{t\leq
T}|\mathbf{u}_{\alpha}(t)|_{2,2},\alpha\in J.
\]

{Then there is a constant }$B_{0}${ such that}%
\[
\tilde{L}_{\alpha}\leq\sqrt{\alpha!}C_{|\alpha|-1}\binom{|\alpha|}{\alpha
}B_{0}^{|\alpha|-1}K^{|\alpha|},|\alpha|\geq2,
\]
{where }$K=1+\sup_{i}\tilde{L}_{\varepsilon_{i}},${ and }%
\[
C_{|\alpha|-1}=\frac{1}{|\alpha|-1}\binom{2(|\alpha|-1)}{|\alpha|-1}%
,|\alpha|\geq2
\]
{are the Catalan numbers (see e.g. \cite{catalan}, \cite{kl} ). Moreover,
there is a number }$q>1${ so that}%
\[
\sup_{r\leq t\leq T}||u(t)||_{_{\mathcal{S}_{-1,-q}}}^{2}\leq\sum_{\alpha
}\frac{(2\mathbf{N)}^{-2q\alpha}\tilde{L}_{\alpha}^{2}}{|\alpha|!}<\infty,
\]
{i.e., the solution }$\mathbf{u}(t)${ belongs to the Kondratiev space of
generalized random functions }$S_{-1}(\mathbb{H}_{p}^{2}\cap\mathbb{H}_{2}%
^{2})=\left(  \mathcal{N}\right)  ^{-1}\otimes(\mathbb{H}_{p}^{2}%
\cap\mathbb{H}_{2}^{2})${. }
\end{proposition}

\begin{proof}
For $|\alpha|\geq1$, $\mathbf{u}_{\alpha}\in C\left(  [0,T],\mathbb{H}_{p}%
^{2}\cap\mathbb{H}_{2}^{2}\right)  ,$ are solutions to the Stokes equations
(\ref{pf4}):%
\begin{align*}
\mathbf{u}_{\alpha}\left(  t\right)   &  =\mathbf{w}_{\alpha}+\int_{r}%
^{t}\mathcal{P}[\partial_{i}\left(  a^{\,ij}\left(  s\right)  \partial
_{j}\mathbf{u}_{\alpha}\left(  s\right)  \right)  +\mathbf{F}_{\alpha}(s)\\
&  +[b^{i}(s)-u_{0}^{i}(s)]\partial_{i}\mathbf{u}_{\alpha}(s)-u_{\alpha}%
^{k}\left(  s\right)  \partial_{k}\mathbf{u}_{0}(s)]ds,
\end{align*}
with $\mathbf{F}_{\alpha}(s)$ defined by (\ref{pf5}).{ Finally, by Proposition
\ref{aprop1} in Appendix II, }%
\[
\tilde{L}_{\alpha}\leq C\sum_{l=2,p}(\int_{0}^{T}|\mathbf{F}_{\alpha
}(s)|_{1,l}^{l}ds)^{1/l}.
\]
Since%
\[
|u_{\alpha-\gamma}^{k}\left(  s\right)  \partial_{k}\mathbf{u}_{\gamma
}(s)|_{1,l}\leq C|\mathbf{u}_{\alpha-\gamma}\left(  s\right)  |_{2,p}%
|\nabla\mathbf{u}_{\gamma}(s)|_{1,l},
\]
it follows that
\begin{align*}
|\mathbf{F}_{\alpha}(s)|_{1,l}  &  \leq C[\sum_{\gamma\leq\alpha
,|\alpha|-1\geq|\gamma|\geq1}\sqrt{\binom{\alpha}{\gamma}}\tilde{L}%
_{\alpha-\gamma}\tilde{L}_{\gamma}\\
&  +\sum_{k}1_{\sigma\neq0}\sqrt{\alpha_{k}}\tilde{L}_{\alpha(k)}%
+1_{|\alpha|=1}|\mathbf{g}(s)|_{1,l}].
\end{align*}
For $|\alpha|\geq2$%
\[
\tilde{L}_{\alpha}\leq C[\sum_{\gamma\leq\alpha,1\leq|\gamma|\leq|\alpha
|-1}\sqrt{\binom{\alpha}{\gamma}}\tilde{L}_{\alpha-\gamma}\tilde{L}_{\gamma
}+1_{\sigma\neq0}\sum_{k}\sqrt{\alpha_{k}}\tilde{L}_{\alpha(k)}].
\]
{ So, there is a constant }$B_{0}${ so that for }$|\alpha|=n\geq2,\hat
{L}_{\alpha}=(\alpha!)^{-1/2}\tilde{L}_{\alpha},\hat{L}_{\varepsilon_{i}%
}=\tilde{L}_{\varepsilon_{i}}${ we have }%
\[
\hat{L}_{\alpha}\leq B_{0}\left(  \sum_{\gamma\leq\alpha,1\leq|\gamma
|\leq|\alpha|-1}\hat{L}_{\alpha-\gamma}\hat{L}_{\gamma}+1_{\sigma\neq0}%
\sum_{k}\hat{L}_{\alpha(k)}1_{\alpha_{k}\neq0}\right)  .
\]
{Denoting }$L_{\alpha}=\hat{L}_{\alpha}${ if }$|\alpha|>1,L_{\alpha}=1+\hat
{L}_{\alpha}${ if }$|\alpha|=1,$ we get%
\[
L_{\alpha}\leq B_{0}\sum_{\gamma\leq\alpha,1\leq|\gamma|\leq|\alpha
|-1}L_{\alpha-\gamma}L_{\gamma}%
\]
{and by \cite{kl} for }$|\alpha|\geq2$%
\begin{align*}
L_{\alpha}  &  \leq C_{|\alpha|-1}B_{0}^{|\alpha|-1}\binom{|\alpha|}{\alpha
}\prod_{i}(1+\tilde{L}_{\varepsilon_{i}})^{\alpha_{i}}\\
&  \leq C_{|\alpha|-1}\binom{|\alpha|}{\alpha}B_{0}^{|\alpha|-1}K^{|\alpha|}%
\end{align*}
with
\[
K=1+C[\tilde{L}_{0}+\sum_{l=2,p}(\int_{0}^{T}||\mathbf{g}(s)||_{1,l}%
^{l}ds)^{1/l}].
\]
{So,}%
\begin{align*}
\tilde{L}_{\alpha}  &  \leq\sqrt{\alpha!}C_{|\alpha|-1}\binom{|\alpha|}%
{\alpha}B_{0}^{|\alpha|-1}K^{|\alpha|},\\
\tilde{L}_{\alpha}^{2}  &  \leq\alpha!C_{|\alpha|-1}^{2}\binom{|\alpha
|}{\alpha}(2\mathbf{N})^{\alpha}B_{0}^{2(|\alpha|-1)}K^{2|\alpha|}%
\end{align*}
{and}%
\[
\frac{r^{\alpha}\tilde{L}_{\alpha}^{2}}{\alpha!}\leq C_{|\alpha|-1}^{2}%
\binom{|\alpha|}{\alpha}(2\mathbf{N}r)^{\alpha}B_{0}^{2(|\alpha|-1)}%
K^{2|\alpha|}.
\]
{Therefore with }$r=(r_{i}),r_{i}=(2i)^{-2q},q>1,$%
\begin{align*}
\sum_{|\alpha|=n}\frac{r^{\alpha}\tilde{L}_{\alpha}^{2}}{\alpha!}  &
=C_{n-1}^{2}B_{0}^{2(n-1)}K^{2n}\sum_{|\alpha|=n}\binom{|\alpha|}{\alpha
}(2\mathbf{N}r)^{\alpha}\\
&  \leq C_{n-1}^{2}B_{0}^{2(n-1)}K^{2n}2^{n}2^{-2qn}\left(  \sum_{i=1}%
^{\infty}\frac{1}{i^{2q-1}}\right)  ^{n}.
\end{align*}
{For large }$n${, the Catalan numbers }%
\[
C_{n-1}\approx\frac{4^{n-1}}{\sqrt{\pi}(n-1)^{3/2}}%
\]
{and there is a number }$q>1${ such that}%
\[
\sum_{n=0}^{\infty}\sum_{|\alpha|=n}\frac{r^{\alpha}\tilde{L}_{\alpha}^{2}%
}{\alpha!}<\infty.
\]
{Therefore, the solution }$\mathbf{u}(t)${ belongs to Kondratiev's space
}$S_{-1}(\mathbb{H}_{p}^{2}\cap\mathbb{H}_{2}^{2}).$
\end{proof}

\begin{remark}
\label{rem3}If Let $\mathbf{u}(t)$ is the solution {(\ref{nse0}) with a
deterministic }$\mathbf{w}\in\mathbb{H}_{p}^{2}\cap\mathbb{H}_{2}^{2},${ and
}$q>1$ is the number in Proposition \ref{cle1}, then the action of
$\mathbf{u}(t)=\sum_{\alpha}\mathbf{u}_{\alpha}(t)\xi_{\alpha}$ can be
extended from $\mathcal{D}(\mathbf{b})$ to%
\[
\mathcal{S}_{1,q}(\mathbf{R)}=\left\{  \eta=\sum_{\alpha}a_{\alpha}\xi
_{\alpha}\in\mathcal{S}_{1}(\mathbf{R}):||\eta||_{\mathcal{S}_{1,q}}%
<\infty\right\}
\]
as%
\[
\left\langle \mathbf{u}(t),\eta\right\rangle =\sum_{\alpha}\mathbf{u}_{\alpha
}(t)a_{\alpha},\eta=\sum_{\alpha}a_{\alpha}\xi_{\alpha}\in\mathcal{S}%
_{1,q}(\mathbf{R)}.
\]

Note that the stochastic exponent $p(z)=p(e_{z})$ (see Remark \ref{rem1})
belongs to $\mathcal{S}_{1,q}(\mathbf{R})$ in (a) provided
\[
|A^{q}e_{z}|_{Y}^{2}=\sum_{k}(z_{k}2^{q}k^{q})^{2}<1.
\]

\end{remark}

{For }$\varepsilon>0${ define a self-adjoint positive operator }%
$D_{\varepsilon}${ on }$H${ such that }$D_{\varepsilon}e_{k}=2^{-\varepsilon
k}e_{k}${ and a sequence of positive numbers }$\kappa_{\varepsilon
,n}=e^{-\varepsilon e^{n}}${. Set }$C_{\varepsilon}=\sum_{n=0}^{\infty}%
\kappa_{\varepsilon,n}D_{\varepsilon}^{\otimes n}.$ It is a second
quantization operator in the Fock space $\mathcal{H}=\sum_{n}\mathcal{H}%
_{n},\mathcal{H}_{n}=\mathbf{H}^{\hat{\otimes}n}$ (see Appendix I, 5.1).
{Clearly, }%
\begin{equation}
C_{\varepsilon}e_{\alpha}=\kappa_{\varepsilon,|\alpha|}D_{\varepsilon
}^{\otimes n}e_{\alpha}=\kappa_{\varepsilon,|\alpha|}\left(  2^{-\varepsilon
\mathbf{N}}\right)  ^{\alpha}e_{\alpha}, \label{eq0}%
\end{equation}
where%
\[
\left(  2^{-\varepsilon\mathbf{N}}\right)  ^{\alpha}=\prod_{k=1}^{\infty
}2^{-\varepsilon k\alpha_{k}}.
\]

\begin{proposition}
{\label{p1}Assume A1-A3 hold and }$\sup_{s,k}|e_{k}(s)|_{Y}<\infty.${ Let
}$\mathbf{u}(t)=\sum_{\alpha}\mathbf{u}_{\alpha}(t)\xi_{\alpha}\in
C\mathcal{D}^{\prime}([0,T],\mathbf{b},\mathbb{H}_{p}^{2}\cap\mathbb{H}%
_{2}^{2})${ be a generalized }$\mathcal{D}$-$\mathbb{H}_{p}^{2}\cap
\mathbb{H}_{2}^{2}${-solution of equation (\ref{nse}) in }$[0,T]${ with a
deterministic }$\mathbf{w}\in\mathbb{H}_{p}^{2}\cap\mathbb{H}_{2}^{2}$ {and}%
\[
\mathbf{u}_{\varepsilon}(t)=C_{\varepsilon}\mathbf{u}(t)=\sum_{n=0}^{\infty
}\kappa_{\varepsilon,n}\sum_{|\alpha|=n}\mathbf{u}_{\alpha}(t)(2^{-\varepsilon
\mathbf{N}})^{\alpha}\xi_{\alpha},
\]
where $C_{\varepsilon}\mathbf{u}$ is rescaling based on the second
quantization operator $C_{\varepsilon}$ (see Appendix I, 5.1).

{Then }$\mathbf{u}_{\varepsilon}(t)${ is }$\mathbb{H}_{2}^{2}${-valued square
integrable process satisfying the equation }%
\begin{align}
\partial_{t}\mathbf{u}_{\varepsilon}\left(  t\right)   &  =\mathcal{P}%
\{\partial_{i}\left(  a^{\,ij}\left(  t\right)  \partial_{j}\mathbf{u}%
_{\varepsilon}\left(  t\right)  \right)  +b^{i}(t)\partial_{i}\mathbf{u}%
_{\varepsilon}(t)\label{for1}\\
&  -C_{\varepsilon}(C_{\varepsilon}{}^{-1}u_{\varepsilon}^{k}\left(  t\right)
)\lozenge(C_{\varepsilon}{}^{-1}\partial_{k}\mathbf{u}_{\varepsilon
}(t))+\mathbf{f}\left(  t\right) \nonumber\\
&  +C_{\varepsilon}[(\sigma^{i}(t)C_{\varepsilon}{}^{-1}\partial_{i}%
\mathbf{u}_{\varepsilon}\left(  t\right)  +\mathbf{g}\left(  t\right)
)\lozenge(C_{\varepsilon}{}^{-1}\dot{W}_{t}^{\varepsilon})]\},\nonumber\\
\mathbf{u}_{\varepsilon}\left(  0\right)   &  =\mathbf{w},\operatorname{div}%
\mathbf{u}_{\varepsilon}=0.\nonumber
\end{align}
{Moreover, }$\mathbf{u}_{\varepsilon}(t)\in\mathcal{S}_{1}(\mathbb{H}_{p}%
^{2}\cap\mathbb{H}_{2}^{2})=\left(  \mathcal{N}\right)  ^{1}\otimes
(\mathbb{H}_{p}^{2}\cap\mathbb{H}_{2}^{2}),t\in\lbrack0,T],\mathbf{u}%
(t)=C_{\varepsilon}{}^{-1}\mathbf{u}_{\varepsilon}(t)${ and }%
\[
\sup_{t\leq T}||\mathbf{u}_{\varepsilon}(t)-\mathbf{u}(t)||_{S_{-1,-q}%
}\rightarrow0{\ }%
\]
{as }$\varepsilon\rightarrow0,${ where }$q${ is a number in Proposition
\ref{cle1}.}
\end{proposition}

\begin{proof}
{Let }$\tilde{L}_{\alpha}=\sup_{t\leq T}|\mathbf{u}_{\alpha}(t)|_{2,p}%
+\sup_{t\leq T}|\mathbf{u}_{\alpha}(t)|_{2,2}$. Since $\mathbf{u}%
_{\varepsilon,\alpha}(t)=\kappa_{\varepsilon,|\alpha|}\left(  2^{-\varepsilon
\mathbf{N}}\right)  ^{\alpha}\mathbf{u}_{\alpha}(t),$
\[
\tilde{L}_{\varepsilon,\alpha}=\sup_{t\leq T}|\mathbf{u}_{\varepsilon,\alpha
}(t)|_{2,p}+\sup_{t\leq T}|\mathbf{u}_{\varepsilon,\alpha}(t)|_{2,2}%
=\kappa_{\varepsilon,|\alpha|}\left(  2^{-\varepsilon\mathbf{N}}\right)
^{\alpha}\tilde{L}_{\alpha}.
\]
{Since for each }$q^{\prime}\geq0,$ there is a constant $C(\varepsilon
,q^{\prime},q)$ independent of $\alpha$ so that
\begin{align*}
&  (|\alpha|!)^{2}\left(  2\mathbf{N}\right)  ^{2(q^{\prime}+q)}%
\kappa_{\varepsilon,|\alpha|}^{2}\left(  2^{-2\varepsilon\mathbf{N}}\right)
^{\alpha}\\
&  \leq\sum_{k}\left(  2k\right)  ^{2(q^{\prime}+q)}2^{-2\varepsilon
k})^{|\alpha|}e^{-2\varepsilon e^{|\alpha|}}(|\alpha|!)^{2}\\
&  \leq C(\varepsilon,q^{\prime},q)<\infty,
\end{align*}
it follows {by Proposition \ref{cle1} that}
\[
||\mathbf{u}_{\varepsilon}(t)||_{\mathcal{S}_{1,q^{\prime}}}^{2}\leq
C(\varepsilon,q^{\prime},q)\sum_{\alpha}\frac{(2\mathbf{N)}^{-2q\alpha}%
\tilde{L}_{\alpha}^{2}}{|\alpha|!}<\infty.
\]
{So, }$\mathbf{u}_{\varepsilon}(t)\in S_{1}(\mathbb{H}_{p}^{2}\cap
\mathbb{H}_{2}^{2}),t\in\lbrack0,T].${ In particular,}%
\[
\mathbf{E}|\mathbf{u}_{\varepsilon}(t)|_{2,2}^{2}=\sum_{\alpha}|\mathbf{u}%
_{\varepsilon,\alpha}(t)|_{2,2}^{2}<\infty,t\in\lbrack0,T].
\]
{Therefore }$\mathbf{u}_{\varepsilon}(t)${ is }$\mathbb{H}_{2}^{2}${-valued
square integrable process and (\ref{for1}) follows by Remark \ref{re4} in
Appendix I, 5.1. Obviously, }$\mathbf{u}(t)=(C_{\varepsilon})^{-1}%
\mathbf{u}_{\varepsilon}(t).${ Since }%
\[
\mathbf{u}_{\varepsilon}(t)-\mathbf{u}(t)=\sum_{\alpha}[1-\kappa
_{\varepsilon,|\alpha|}\left(  2^{-\varepsilon\mathbf{N}}\right)  ^{\alpha
}]\mathbf{u}_{\alpha}(t,x)\xi_{\alpha},
\]
it follows{ that}%
\[
||\mathbf{u}_{\varepsilon}(t)-\mathbf{u}(t)||_{\mathcal{S}_{-1,-q}}^{2}%
\leq\sum_{\alpha}|1-\kappa_{\varepsilon,|\alpha|}\left(  2^{-\varepsilon
\mathbf{N}}\right)  ^{\alpha}|^{2}\frac{(2\mathbf{N)}^{-2q\alpha}\tilde
{L}_{\alpha}^{2}}{|\alpha|!}\rightarrow0
\]
{as }$\varepsilon\rightarrow0${ by Lebesgue's dominated convergence theorem
uniformly in }$t${.}
\end{proof}

\begin{remark}
\label{rem4}The solution $\mathbf{u}$ in Proposition \ref{cle1} depends on a
fixed uniformly bounded basis $\mathbf{b}$ in $\mathbf{H}$. It belongs to the
Kondratiev space
\[
S_{-1}(\mathbb{H}_{p}^{2}\cap\mathbb{H}_{2}^{2})=\left(  \mathcal{N}\right)
^{-1}\otimes(\mathbb{H}_{p}^{2}\cap\mathbb{H}_{2}^{2})
\]
constructed using a Gelfand triple $\mathcal{N}\subseteq$ $\mathbf{H}%
=L_{2}\left(  [0,T],Y\right)  \subseteq\mathcal{N}^{\prime}$ which depends on
$\mathbf{b}$.

In \cite{gks}, a class $\mathcal{G}^{-1}\subseteq\left(  \mathcal{N}\right)
^{-1}$ of regular generalized functions was introduced that does not depend on
a fixed Gelfand triple or a basis in $\mathbf{H}$. Unfortunately, the
estimates in Proposition \ref{cle1} (because of the factor $\mathbf{N}%
^{-q\alpha}=\Pi_{k}k^{-q\alpha_{k}}$) do not imply that $\mathbf{u}$ is a
regular generalized function of class $\mathcal{G}^{-1}$. Also, the space
\[
\mathcal{S}_{-1,-q}(\mathbb{H}_{p}^{2}\cap\mathbb{H}_{2}^{2})=\left\{
\mathbf{\eta}\in\mathcal{S}_{-1}(\mathbb{H}_{p}^{2}\cap\mathbb{H}_{2}%
^{2}):||\mathbf{\eta}||_{\mathcal{S}_{-1,-q}}<\infty\right\}  ,
\]
with some $q>1,$ to which the solution in Proposition \ref{cle1} belongs,
cannot be embedded into any space with weights depending only on $|\alpha|$
(for example, into the spaces, like $\mathcal{G}^{-1,-q}$ in \cite{gks}$,$
that do not dependent on a fixed \ basis or Gelfand triple in $\mathbf{H}$):
$\inf_{|\alpha|=n}\left(  \mathbf{N}\right)  ^{-q\alpha}=0.$
\end{remark}

\section{Markov property and independence of basis}

In this Section we will show that a generalized $\mathcal{D}$-$\mathbb{H}%
_{p}^{2}\cap\mathbb{H}_{2}^{2}$-solution of equation (\ref{nse0}) has the
following properties: it is adapted with respect to the {filtration
}$(\mathcal{F}_{t}^{W})$ {generated by the Wiener process }$W_{t};${ it is
independent of the choice of the basis }$\mathbf{b}${, and it is a generalized
Markov process.}

\subsection{Equivalent characterization of $\mathcal{D}$-generalized
processes}

A more convenient characterization of $\mathcal{D}$-solution to\ (\ref{2})
(see Definition \ref{defs}) is based on another (equivalent) description of
$\mathcal{D}(\mathbf{b).}$ It allows to introduce the notion of an adapted
solution and extend it from $\mathcal{D}(\mathbf{b)}$ to a space of test
functions that is independent of $\mathbf{b}\in\mathcal{B}$.

\subsubsection{Equivalent description of test function space}

Often it is convenient to use the exponents $p(z),z\in\mathcal{Z},$ defined in
Remark \ref{rem1} to describe the test function space $\mathcal{D}%
(\mathbf{b)}$.

{To each multi-index }$\alpha${ of length }$n${ we relate a set }$K_{\alpha}${
whose elements are positive integers }$k_{i},i=1,\ldots,n,${ such that each
}$k${ is represented there by }$\alpha_{k}${-copies. An ordered }$n${-tuple
}$K_{\alpha}=\{k_{1},\ldots,k_{n}\}$ with $k_{1}\leq k_{2}\leq\ldots\leq
k_{n}${ characterizes the locations and the values of the non-zero components
of }$\alpha${. For example, }$k_{1}${ is the index of the first non-zero
element of }$\alpha,${ followed by }$\max\left(  0,\alpha_{k_{1}}-1\right)  ${
of entries with the same value (see \cite{MRWC}). }

{For an orthonormal basis }$\{e_{k},k\geq1\}${\ in }$L_{2}\left(
[0,T],Y\right)  ${ and }$\alpha\in I${ with }$K_{\alpha}=\left\{  k_{1}%
,\ldots,k_{n}\right\}  ${, we denote}%
\[
E_{\alpha}=\sum_{\sigma\in G^{n}}e_{k_{\sigma(1)}}\otimes\ldots\otimes
e_{k_{\sigma(n)}},\alpha\in J,
\]
{where }$G^{n}${ is a permutation group of }$\left\{  1,\ldots,n\right\}  ${.
The set}%
\begin{equation}
\left\{  e_{\alpha}=\frac{E_{\alpha}}{\sqrt{\alpha!|\alpha|!}},\alpha\in
J\right\}  \label{f0}%
\end{equation}
{ is a CONS for the symmetric part }$\mathbf{H}^{\hat{\otimes}n}$ {of
}$\mathbf{H}^{\otimes n}.$

{For }$|\alpha|=n,${ }%
\begin{equation}
\xi_{\alpha}=\sqrt{|\alpha|!}W(e_{\alpha}), \label{Eq:XiAlpha}%
\end{equation}
{where }%
\begin{equation}
W(e_{\alpha})=\int_{0}^{T}\int_{0}^{s_{n}}\ldots\int_{0}^{s_{2}}e_{\alpha
}(s_{1},\ldots,s_{n})dW_{s_{1}}\ldots dW_{s_{n}}. \label{Eq:WeAlpha}%
\end{equation}
If{ }$\alpha=\varepsilon_{k}${, then }$W(e_{\varepsilon_{k}})=W\left(
e_{k}\right)  =\int_{0}^{T}e_{k}(t)dW_{t}$.

According to (\ref{f4}),
\begin{align*}
p(z)  &  =p(e_{z})=p(z,\mathbf{b})=\sum_{n=0}^{\infty}\sum_{|\alpha
|=n}z^{\alpha}\sqrt{\frac{|\alpha|!}{\alpha!}}W(e_{\alpha})\\
&  =\sum_{n=0}^{\infty}\sum_{|\alpha|=n}\frac{z^{\alpha}}{\sqrt{\alpha!}}%
\xi_{\alpha}.
\end{align*}
{Denote }%
\begin{align}
p_{n}(z)  &  =p_{n}(e_{z})=\sum_{|\alpha|=n}\frac{z^{\alpha}}{\sqrt{\alpha!}%
}\xi_{\alpha}\nonumber\\
&  =\int_{0}^{T}\int_{0}^{s_{n-1}}\ldots\int_{0}^{s_{2}}e_{z}(s_{1})\ldots
e_{z}(s_{n})dW_{s_{1}}\ldots dW_{s_{n}},n\geq2,\label{f5}\\
p_{1}(z)  &  =p_{1}(e_{z})=\int_{0}^{T}e_{z}(s_{1})dW_{s_{1}},p_{0}%
(z)=p_{0}(e_{z})=1.\nonumber
\end{align}

\begin{lemma}
{\label{TeqD}For }$\mathbf{b\in}\mathcal{B}$, {let }$\mathcal{V}${
=}$\mathcal{V}(\mathbf{b})${ be the linear space of random variables that
consists of all finite linear combinations of }$p_{n}(z),z\in\mathcal{Z}%
,n\geq0${. Then }$\mathcal{V}(\mathbf{b})=\mathcal{D}(\mathbf{b})$ (in
particular, $\xi_{\alpha}\in\mathcal{V}(\mathbf{b})$)$.$
\end{lemma}

\begin{proof}
{Obviously, }$\mathcal{V}=\mathcal{V}(\mathbf{b})\subseteq\mathcal{D}%
=\mathcal{D}(\mathbf{b}).${ For }$\alpha\in\mathcal{I}${, denote }%
$\kappa(\alpha)=\max\left\{  k:\alpha_{k}\neq0\right\}  ${. Fix }$N,n${ and
}$\alpha=\left(  \alpha_{k}\right)  \in\mathcal{I}${ such that }%
$|\alpha|=N,\kappa(\alpha)=n${. Consider a finite dimensional Hilbert space}%
\[
G=\left\{  \sum_{|\alpha|=N,\kappa(\alpha)\leq n}v_{\alpha}\xi_{\alpha
}:v_{\alpha}\in\mathbf{R}\right\}  .
\]
{with inner product }%
\[
(\sum_{a}v_{\alpha}\xi_{\alpha},\sum_{a}v_{\alpha}^{\prime}\xi_{\alpha}%
)_{G}=\sum_{\alpha}v_{\alpha}v_{\alpha}^{\prime}.
\]
Let $\tilde{G}$ be {a vector subspace of }$G${ generated by }$p_{N}%
(z),z=(z_{1},\ldots,z_{n},0,\ldots)\in\mathcal{Z}.${ It is enough to show that
}$\tilde{G}=G${. Indeed, the subspace }$\tilde{G}${ is }finite-dimensional{
and obviously closed. Assume there is a vector }$\sum_{\alpha}v_{\alpha}%
\xi_{\alpha}\in G${ which is orthogonal to }$\tilde{G}.${ So, for all
}$z=(z_{1},\ldots,z_{n},0,\ldots)\in\mathcal{Z},$%
\[
(\sum_{\alpha}v_{\alpha}\xi_{\alpha},p_{N}(z))_{G}=\sum_{\alpha}v_{\alpha
}\frac{z^{\alpha}}{\sqrt{\alpha!}}=0
\]
{which implies that all }$v_{\alpha}=0.${ Therefore }$\tilde{G}=G.${This
completes the proof.}
\end{proof}

{Due to Lemma \ref{TeqD}, we can characterize convergence in }$\mathcal{D}%
^{\prime}=\mathcal{D}^{\prime}(\mathbf{b)}${ by test functions of the form
}$p_{m}\left(  z\right)  ${. Indeed, for }$z\in\mathcal{Z}$, $v\in\mathcal{D}%
$, and $m\geq0${, we have}%
\begin{equation}
\left\langle p_{m}\left(  z\right)  ,v\right\rangle =\sum_{|\alpha
|=m}v_{\alpha}\frac{z^{\alpha}}{\sqrt{\alpha!}}. \label{ff1}%
\end{equation}
{Therefore we have the following necessary and sufficient condition:}

\begin{corollary}
{\label{FF1}\bigskip A sequence }$v^{n}\rightarrow v${ in }$\mathcal{D}%
^{\prime}${ if and only if for all }$z\in\mathcal{Z}${ and }$\ ${all }$m\geq
0${ }%
\[
\left\langle p_{m}\left(  z\right)  ,v^{n}\right\rangle \rightarrow
\left\langle p_{m}\left(  z\right)  ,v\right\rangle .
\]

\end{corollary}

\subsubsection{Action of a Skorokhod integral on $p_{M}(z)$}

{Consider }$v(t)=\sum_{\alpha}v_{\alpha}(t)\xi_{\alpha}\in\mathcal{D}^{\prime
}\left(  \mathbf{b};[0,T],Y\right)  ${ such that for all }$\alpha,k,$
$\int_{0}^{T}|\left(  v_{\alpha}(s),e_{k}(s)\right)  _{Y}|ds<\infty.$ Recall
that the {Skorokhod integral assigns to\ such }$v${ a generalized random
}process
\begin{align*}
\delta_{t}(v)  &  =\int_{0}^{t}v(s)dW_{s}=\delta\left(  v1_{[0,t]}\right)
=\sum_{\alpha}\delta_{t}(v)_{\alpha}\xi_{\alpha},0\leq t\leq T,\\
\delta_{t}(v)_{\alpha}  &  =\sum_{k}\sqrt{\alpha_{k}}\int_{0}^{t}\left(
v_{\alpha(k)}(s),e_{k}(s)\right)  _{Y}ds.
\end{align*}

\begin{remark}
{\label{re2}For }$p_{M}(z)\in\mathcal{D}(\mathbf{b}),z\in Z,M\geq1,$ it is
easy to show that%
\[
\left\langle p_{M}\left(  z\right)  ,\delta_{t}(v)\right\rangle =\int_{0}%
^{t}(\left\langle p_{M-1}\left(  z\right)  ,v(s)\right\rangle ,e_{z}%
(s))_{Y}~ds.
\]

\end{remark}

\subsubsection{Action of a Wick product on $p_{M}(z)$}

Recall that {for a Hilbert space }$E${ and arbitrary }$v=\sum_{\alpha
}v_{\alpha}\xi_{\alpha}${ and }$u=\sum_{\alpha}u_{\alpha}\xi_{\alpha}${ in
}$\mathcal{D}^{\prime}(\mathbf{b},E)${, we define }%
\[
v\lozenge u=\sum_{\alpha}\sum_{\beta\leq\alpha}(u_{\beta},v_{\alpha-\beta
})_{E}\sqrt{\frac{\alpha!}{\beta!(\alpha-\beta)!}}\xi_{\alpha}\in
\mathcal{D}^{\prime}\mathcal{(}\mathbf{b},\mathbf{R}).
\]

{In particular, }%
\[
\xi_{\alpha}\lozenge\xi_{\beta}=\xi_{\alpha+\beta}\sqrt{\frac{(\alpha+\beta
)!}{\beta!\alpha!}}.
\]

{The following statement holds.}

\begin{lemma}
{\label{le2}For a Hilbert space }$E,${ arbitrary elements }$v=\sum_{\alpha
}v_{\alpha}\xi_{\alpha}${ and }$u=\sum_{\alpha}u_{\alpha}\xi_{\alpha}${ from
}$\mathcal{D}^{\prime}(\mathbf{b},E),${ and}$~z\in\mathcal{Z},M\geq0,$%
\[
\left\langle p_{M}(z),v\lozenge u\right\rangle =\sum_{K+L=M}\left\langle
p_{K}(z),v\right\rangle \left\langle p_{L}(z),u\right\rangle ;
\]
{In particular, }%
\[
\left\langle 1,v\lozenge u\right\rangle =\left(  \left\langle 1,v\right\rangle
,\left\langle 1,u\right\rangle \right)  _{E}%
\]
{(the generalized expected value of }$v\lozenge u${ is the product of expected
values})$.$
\end{lemma}

\begin{proof}
{According to (\ref{ff1}),}%
\begin{align*}
\left\langle p_{M}(z),v\lozenge u\right\rangle  &  =\sum_{|\alpha|=M}%
\sum_{\beta\leq\alpha}(\frac{z^{\beta}}{\sqrt{\beta!}}u_{\beta},\frac
{z^{\alpha-\beta}}{\sqrt{(\alpha-\beta)!}}v_{\alpha-\beta})_{E}\\
&  =\left(  \sum_{K+L=M}\sum_{|\beta|=K}\frac{z^{\beta}}{\sqrt{\beta!}%
}u_{\beta},\sum_{|\gamma|=L}\frac{z^{\gamma}}{\sqrt{\gamma!}}u_{\gamma
}\right)  _{E}\\
&  =\left(  \left\langle p(z),v\right\rangle ,\left\langle p(z),u\right\rangle
\right)  _{E}.
\end{align*}

\end{proof}

\subsubsection{An equivalent characterization of the solution}

{Now, we will characterize the solution of equation (\ref{nse}) by its action
on test functions }$p_{M}(z),z\in\mathcal{Z},M\geq0.${\ The following
statement holds.}

\begin{remark}
{\label{c1}Assume A1-A3 hold, }$\mathbf{w}=\sum_{\alpha}\mathbf{w}_{\alpha}%
\xi_{\alpha}\in\mathcal{D}^{\prime}(\mathbf{b},\mathbb{H}_{p}^{2}%
\cap\mathbb{H}_{2}^{2}),\operatorname{div}\mathbf{w}=0${ and }$\mathbf{u}%
(t)=\sum_{\alpha}\mathbf{u}_{\alpha}(t)\xi_{\alpha}\in C\mathcal{D}^{\prime
}(\mathbf{b};[r,T],\mathbb{H}_{p}^{2}\cap\mathbb{H}_{2}^{2}).$ {Then
}$\mathbf{u}(t)${ is }$\mathcal{D}${-}$\mathbb{H}_{p}^{2}\cap\mathbb{H}%
_{2}^{2}${ solution of (\ref{nse}) in }$[r,T]${ if and only if for all }%
$z\in\mathcal{Z}${ and }$M\geq0,${ }%
\[
\mathbf{u}^{M,z}(t)=\left\langle \mathbf{u}(t),p_{M}(z)\right\rangle
=\sum_{|\alpha|=M}\frac{\mathbf{u}_{\alpha}(t)z^{\alpha}}{\sqrt{\alpha!}}%
\]
{is continuous in }$\mathbb{H}_{p}^{2}\cap\mathbb{H}_{2}^{2}$ and the
following equality holds in $\mathbb{L}_{l},l=2,p,${ }%
\begin{align}
\mathbf{u}^{M,z}\left(  t\right)   &  =\mathbf{w}^{M,z}+\int_{r}%
^{t}\mathcal{P}[\partial_{i}\left(  a^{\,ij}\left(  s\right)  \partial
_{j}\mathbf{u}^{M,z}\left(  s\right)  \right) \label{pf2}\\
&  +b^{i}(s)\partial_{i}\mathbf{u}^{M,z}(s)-\sum_{K+L=M}u^{k,K,,z}\left(
s\right)  \partial_{k}\mathbf{u}^{L,z}(s)\nonumber\\
&  +1_{M=0}\mathbf{f}\left(  s\right)  +\left(  \sigma^{i}(s),e_{z}(s)\right)
_{Y}\partial_{i}\mathbf{u}^{M-1,z}(s)\nonumber\\
&  +1_{M=1}(\mathbf{g}\left(  s\right)  ,e_{z}(s))_{Y}]ds,\nonumber
\end{align}
{where }$M\geq0,${ }$\mathbf{w}^{M,z}(x)=\left\langle \mathbf{w}%
(x),p_{M}(z)\right\rangle =\sum_{|\alpha|=M}\mathbf{w}_{\alpha}(x)z^{\alpha
}/\sqrt{\alpha!}${, and }$\mathbf{u}^{-1,z}(t,x)=0.$ {If }$M\geq1,${ equation
(\ref{pf2}) is Stokes equation; if }$M=0,${ it is Navier-Stokes equation.
Indeed, we obtain (\ref{pf2}) by multiplying both sides of (\ref{pr1}) by
}$z^{\alpha}/\sqrt{\alpha!}${ and adding. }
\end{remark}

\subsection{{Adapted and }independent of basis {generalized processes}}

Let $L_{\infty}\left(  \left[  0,T\right]  ,Y\right)  $ be the space of
measurable $Y$-valued bounded functions on $[0,T]$. For $h\in L_{\infty
}\left(  \left[  0,T\right]  ,Y\right)  ,M\geq0,$ we denote%
\[
p_{M,t}(h)=\int_{0}^{t}\int_{0}^{s_{M}}\ldots\int_{0}^{s_{2}}h(s_{1})\ldots
h(s_{M})dW_{s_{1}}\ldots dW_{s_{M}},0\leq t\leq T.
\]
By (\ref{f5}), $p_{M}(z)=p_{M}(e_{z})=p_{M,T}(e_{z}),z\in\mathcal{Z}.$

\begin{lemma}
\label{re1 copy(1)}(i) If $\{m_{k},k\geq1\}$ is a CONS in $L_{2}\left(
0,T\right)  ,$ and $\{\ell_{k},k\geq1\}$ is a CONS in $Y$,$h\in L_{\infty
}\left(  \left[  0,T\right]  ,Y\right)  $, then for each $n,n^{\prime}\geq1,$
there is $z\in\mathcal{Z}$ such that%
\[
h_{n,n^{\prime}}(t)=\sum_{i=1}^{n^{\prime}}\sum_{k=1}^{n}\int_{0}%
^{T}(h(s),\ell_{k})m_{i}(s)ds\ell_{k}m_{i}(t)=e_{z}(t),
\]
$0\leq t\leq T$. Obviously, $p_{M,T}(h_{n,n^{\prime}})\in\mathcal{D}%
(\mathbf{b})$, $\mathbf{b}=\left\{  e_{k}=m_{i_{k}}\ell_{j_{k}},k\geq
1\right\}  ,$%
\begin{align*}
h_{n,n^{\prime}}  &  \rightarrow h\text{ in }L_{2}([0,T],Y),\\
p_{M,T}(h_{n,n^{\prime}})  &  \rightarrow p_{M,T}(h)\text{ in }L_{2}%
(\Omega,\mathbf{P}),
\end{align*}
as $n,n^{\prime}\rightarrow\infty.$

(ii) Assume $\left(  m_{i}\right)  $ is trigonometric basis or unconditional
$L_{p}\left(  [0,T]\right)  $-basis (for example, Haar basis, see \cite{ym}),
$h\in L_{\infty}([0,T],Y)$. Then there is a sequence $z(n)\in\mathcal{Z}$ such
that $e_{z(n)}\rightarrow h$ in $L_{p}\left(  [0,T],Y\right)  $ for all
$p\geq2$, as $n\rightarrow\infty$.
\end{lemma}

\begin{proof}
We prove the second part of the statement. Let
\[
h_{n}(s)=\sum_{k=1}^{n}(h(s),\ell_{k})_{Y}\ell_{k},n\geq1.
\]
Then $|h_{n}(s)|_{Y}\leq|h(s)|_{Y}$ and for all $p\geq2,$%
\[
\int_{0}^{T}|h_{n}(s)-h(s)|_{Y}^{p}ds\rightarrow0
\]
as $n\rightarrow\infty$. If $\left(  m_{i}\right)  $ is trigonometric basis or
unconditional $L_{p}\left(  [0,T]\right)  $-basis (for example, Haar basis),
then for each $n$%
\[
\int_{0}^{T}|h_{n,n^{\prime}}(s)-h_{n}(s)|_{Y}^{p}ds\rightarrow0
\]
as $n^{\prime}\rightarrow$ $\infty$. So, there is a subsequence $l_{n}$ such
that%
\[
\int_{0}^{T}|h_{n,l_{n}}(s)-h(s)|_{Y}^{p}ds\rightarrow0
\]
as $n\rightarrow\infty$, and (ii) follows according to part (i) of this remark.
\end{proof}

Let $\mathcal{T}$ $\ $be the space of all linear combinations of
$p_{M,T}(h),h\in L_{\infty}([0,T],Y)$, $M\geq0.$ Obviously, $\cup
_{\mathbf{b}\in\mathcal{B}}\mathcal{D}(\mathbf{b})\subseteq\mathcal{T},$ and
$\mathcal{T}$ does not depend on any particular $\mathbf{b}\in\mathcal{B}$.

We say that $h_{n}\rightarrow h$ in $L_{\infty}([0,T],Y)$ if $h_{n}\rightarrow
h$ in $L_{p}([0,T],Y)$ for all $p\geq2$.

\begin{definition}
Given a Banach space $E$, let $\mathcal{T}^{\prime}=\mathcal{T}^{\prime}(E)$
be the space of all linear $E$-valued functions $v$ on $\mathcal{T}$ such that
$v(p_{M,T}(h_{n}))\rightarrow v(p_{M,T}(h))$ for all \thinspace$M\geq0$ if
$h_{n}\rightarrow h$ in $L_{\infty}([0,T],Y)$. We say $v\in\mathcal{T}%
^{\prime}(E)$ is a generalized $E$-valued r.v. $\ $
\end{definition}

Denote $\mathcal{T}^{\prime}(\mathbf{b})=\mathcal{T}^{\prime}(\mathbf{b,}E)$
the set of all $v\in\mathcal{D}^{\prime}(\mathbf{b},E)$ such that for each
$h\in L_{\infty}([0,T],Y)$ and any sequence $e_{z(n)}\rightarrow h$ in
$L_{\infty}([0,T],Y),$ the limit $\lim_{n\rightarrow\infty}\left\langle
p_{M,T}\left(  z(n)\right)  ,v\right\rangle $ exists in $E$ for all $M\geq0,$
and does not depend on a particular sequence $z(n)$ such that $e_{z(n)}%
\rightarrow h$ in $L_{\infty}([0,T],Y).$ We define%
\[
\left\langle p_{M,T}\left(  h\right)  ,v\right\rangle =\lim_{n\rightarrow
\infty}\left\langle p_{M}\left(  z(n)\right)  ,v\right\rangle .
\]

\begin{lemma}
\label{les1}Let $\mathbf{b}\in\mathcal{B}$. Then

(a) $\mathcal{T}^{\prime}(\mathbf{b)\subseteq}\mathcal{T}^{\prime}(E).$

(b) For $v\in\mathcal{T}^{\prime}(E)$, there are $a_{\alpha}\in E$ such that
the restriction%
\[
v|_{\mathcal{D}(\mathbf{b)}}=\sum_{\alpha}a_{\alpha}\xi_{\alpha}%
(\mathbf{b)}\in\mathcal{T}^{\prime}(\mathbf{b},E)\mathbf{.}%
\]

\end{lemma}

\begin{proof}
(a) Let $v\in\mathcal{T}^{\prime}(\mathbf{b)}$, $h_{n}\rightarrow h$ in
$L_{\infty}([0,T],Y)$. For each $n$ there is a sequence $z_{k}(n)$ such that
$e_{z_{k}(n)}\rightarrow h_{n}$ in $L_{\infty}[0,T]$ and $v(z_{k}%
(n))\rightarrow v(h_{n})$ as $k\rightarrow\infty$. Therefore there is a
subsequence $z_{k_{n}}(n)$ such that%
\[
|e_{z_{k_{n}}(n)}-h_{n}|_{L_{n}([0,T],Y)}+|v(e_{z_{k_{n}}(n)})-v(h_{n}%
)|\leq1/n.
\]
Then for each $p\geq2,n\geq p,$%
\begin{align*}
&  |e_{z_{k_{n}}(n)}-h|_{L_{p}([0,T],Y)}\\
&  \leq|e_{z_{k_{n}}(n)}-h_{n}|_{L_{p}([0,T],Y)}+|h_{n}-h|_{L_{p}([0,T],Y)}\\
&  \leq T^{\frac{1}{p}-\frac{1}{n}}|e_{z_{k_{n}}(n)}-h_{n}|_{L_{n}%
([0,T],Y)}+|h_{n}-h|_{L_{p}([0,T],Y)}%
\end{align*}
So, $e_{z_{k_{n}}(n)}\rightarrow h$ in $L_{\infty}\left(  [0,T],Y\right)  $
and%
\begin{align*}
&  |v(h_{n})-v(h)|\\
&  \leq|v(h_{n})-v(e_{z_{k_{n}(n)}})|+|v(e_{z_{k_{n}(n)}})-v(h)|\rightarrow0
\end{align*}
as $n\rightarrow\infty$.

(b) Let $v\in\mathcal{T}^{\prime}(E),\mathbf{b}=\left\{  e_{n}\right\}
,\xi_{\alpha}=\xi_{\alpha}(\mathbf{b}),N\geq1,$%
\[
p_{N}(z)=p_{N}(z,\mathbf{b})=p_{N}(e_{z})=\sum_{|\alpha|=N}\frac{z^{\alpha}%
}{\sqrt{\alpha!}}\xi_{\alpha}.
\]
By Lemma \ref{TeqD}, $\xi_{\alpha}\in\mathcal{T}$ and by linearity%
\[
v(p_{N}(z))=\sum_{|\alpha|=N}\frac{z^{\alpha}}{\sqrt{\alpha!}}v(\xi_{\alpha
}).
\]
Denoting $v_{\alpha}=v(\xi_{\alpha})$ set $\bar{v}=\sum_{\alpha}v_{a}%
\xi_{\alpha}\in\mathcal{D}^{\prime}(\mathbf{b}).$Obviously,%
\[
v(p_{N}(z))=\sum_{|\alpha|=N}\frac{z^{\alpha}}{\sqrt{\alpha!}}v(\xi_{\alpha
})=\bar{v}(p_{N}(z))
\]
and (b) holds.
\end{proof}

Let $\mathcal{T}^{\prime}(\mathbf{b};[r,T])=\mathcal{T}^{\prime}%
(\mathbf{b};[r,T],E)$ be the space of all $v\in\mathcal{D}^{\prime}%
(\mathbf{b};[r,T],E)$ such that $v(t)\in\mathcal{T}^{\prime}(\mathbf{b,}%
E),r\leq t\leq T.$

\begin{definition}
Given a Banach space $E$, let $\mathcal{T}^{\prime}=\mathcal{T}^{\prime
}([r,T],E)$ be the space of all $\mathcal{T}^{\prime}(E)$-valued functions
$v(t)$ on $[r,T]$. We say $v\in\mathcal{T}^{\prime}([r,T],E)$ is a generalized
$E$-valued stochastic process. We denote $C\mathcal{T}^{\prime}([r,T],E)$ the
set of all continuous $u\in\mathcal{T}^{\prime}([r,T],E)$.
\end{definition}

The following obvious consequence of Lemma \ref{les1} holds.

\begin{corollary}
$\ $\label{cors}Let $\mathbf{b}\in\mathcal{B},r<T$. Then

(a) $\mathcal{T}^{\prime}(\mathbf{b;}[r,T],E\mathbf{)\subseteq}\mathcal{T}%
^{\prime}([r,T],E).$

(b) For $v\in\mathcal{T}^{\prime}([r,T],E)$, there are $E$-valued functions
$a_{\alpha}(t),0\leq t\leq T,$ such that the restriction%
\[
v|_{\mathcal{D}(\mathbf{b)}}=\sum_{\alpha}a_{\alpha}(t)\xi_{\alpha
}(\mathbf{b)}\in\mathcal{T}^{\prime}(\mathbf{b;}[r,T],E)\mathbf{.}%
\]

\end{corollary}

Now, we introduce the notion of an adapted generalized process. { }

\begin{definition}
(a) We say $v\in\mathcal{T}^{\prime}(E)$ is $\mathcal{F}_{t_{0}}^{W}%
$-measurable if
\[
\left\langle p_{M,T}(h),v\right\rangle =\left\langle p_{M,t_{0}}%
(h),v\right\rangle
\]
for all $h\in L_{\infty}\left(  \left[  0,T\right]  ,Y\right)  ,${ }$M\geq0.$

(b) {We say }$v\in\mathcal{T}^{\prime}([0,T],E)${ is }$\mathbb{F}^{W}%
${-adapted if }$v(t)$ is $\mathcal{F}_{t}^{W}$-measurable {for each }$t$.
\end{definition}

\begin{example}
{Let }$W_{t}${ be a cylindrical Wiener process in a Hilbert space }$Y${ and
}$\dot{W}_{t}=${ }$\frac{d}{dt}W_{t}.${ Then (see Example \ref{ex:CWP} as
well) }$W_{t}${ and }$\dot{W}_{t}${ are generalized }$Y${-valued }adapted
{stochastic processes. For any }$h\in L_{\infty}(\left[  0,T\right]  ,Y),$%
\begin{align*}
\left\langle W_{t},p_{M,T}(h)\right\rangle  &  =\int_{0}^{t}%
h(s)ds=\left\langle W_{t},p_{M,t}(h)\right\rangle ,\\
\left\langle \dot{W}_{t},p_{M,T}(h)\right\rangle  &  =h(t)=\left\langle
\dot{W}_{t},p_{M,t}(h)\right\rangle
\end{align*}
{if }$M=1${, and }$\left\langle W_{t},p_{s}^{M}(h)\right\rangle =\left\langle
\dot{W}_{t},p_{s}^{M}(h)\right\rangle =0,0\leq s\leq T,${ otherwise.}
\end{example}

\subsection{Independence of basis and Markov property of the solution}

Remark \ref{c1} {suggests the following definition of a generalized solution
to }(\ref{nse0}){.}

\begin{definition}
{Given }$\mathbf{w}\in\mathcal{T}^{\prime}(\mathbb{H}_{p}^{2}\cap
\mathbb{H}_{2}^{2}),\operatorname{div}\mathbf{w}=0,T\geq r,${ a generalized
process }$\mathbf{u}\in C\mathcal{T}^{\prime}([r,T],\mathbb{H}_{p}^{2}%
\cap\mathbb{H}_{2}^{2})${ is called }$\mathbb{H}_{p}^{2}\cap\mathbb{H}_{2}%
^{2}${-solution of equation (\ref{nse}) in }$[r,T]${, if for each }$h\in
L_{\infty}([0,T],Y),M\geq0,${ the function }$\mathbf{u}^{M,h}%
(t,x)=\left\langle p_{M,T}(h),\mathbf{u}(t,x)\right\rangle ${ is an
}$\mathbb{H}_{p}^{2}\cap\mathbb{H}_{2}^{2}$-valued continuous functions
satisfying {in }$\mathbb{L}_{l},l=2,p,$ {Stokes (}$M\geq1)${ or Navier-Stokes
(}$M=0)${ equation}%
\begin{align}
\mathbf{u}^{M,h}\left(  t\right)   &  =\mathbf{w}^{M,h}+\int_{r}%
^{t}\mathcal{P}[\partial_{i}\left(  a^{\,ij}\left(  s\right)  \partial
_{j}\mathbf{u}^{M,h}\left(  s\right)  \right) \label{pf3}\\
&  +b^{i}(s)\partial_{i}\mathbf{u}^{M,h}(s)-\sum_{K+L=M}u^{k,K,,h}\left(
s\right)  \partial_{k}\mathbf{u}^{L,h}(s)\nonumber\\
&  +1_{M=0}\mathbf{f}\left(  s\right)  +\left(  \sigma^{i}(s),h(s)\right)
_{Y}\partial_{i}\mathbf{u}^{M-1,h}(s)\nonumber\\
&  +1_{M=1}(\mathbf{g}\left(  s\right)  ,h(s))_{Y}]ds,\nonumber
\end{align}
{where }$M\geq0${, }$\mathbf{w}^{M,h}(x)=\left\langle p_{M,T}(h),\mathbf{w}%
(x)\right\rangle ${, and }$\mathbf{u}^{-1,h}(t,x)=0.$
\end{definition}

{Obviously, a generalized solution is a }$\mathcal{D}${-solution. }Now we are
in a position to prove the main result.

\begin{theorem}
{\label{main}Assume that A1-A3 hold, }$\mathbf{w}\in\mathcal{T}^{\prime
}(\mathbb{H}_{p}^{2}\cap\mathbb{H}_{2}^{2})${. Then for each }$r<T<T_{1}${
there is a unique }$\mathbb{H}_{p}^{2}\cap\mathbb{H}_{2}^{2}${-solution
}$\mathbf{u}\in C\mathcal{T}^{\prime}([r,T],\mathbb{H}_{p}^{2}\cap
\mathbb{H}_{2}^{2})${\ of (\ref{nse0}) in }$[0,T].${ Moreover, if
}${\mathbf{\ w}}$ is $\mathcal{F}_{r}^{W}$-measurable, then {the solution is
}$\mathbb{F}^{W}=\left(  \mathcal{F}_{t}^{W}\right)  _{t\geq r}${-adapted and
it extends all }$\mathcal{D}(\mathbf{b})${-solutions, }$\mathbf{b}%
\in\mathcal{B}${.}
\end{theorem}

Before proceeding with the proof of Theorem \ref{main} we shall prove the
following auxiliary statement.

\begin{lemma}
\label{a1}Let {A1-A3 hold, }$\mathbf{w}\in\mathcal{T}^{\prime}(\mathbb{H}%
_{p}^{2}\cap\mathbb{H}_{2}^{2}),\operatorname{div}\mathbf{w}=0,h\in L_{\infty
}\left(  [0,T],Y\right)  ${. Then the infinite system (\ref{pf3}) with
}${\mathbf{w}}^{M,h}=\left\langle \mathbf{w},p_{M,T}(h)\right\rangle $ {has a
unique solution }$\mathbf{v}^{M,h}\in C\left(  [r,T],\mathbb{H}_{p}^{2}%
\cap\mathbb{H}_{2}^{2}\right)  ,M\geq0.$ Moreover, for each $N\geq1$ there is
a constant $C$ independent of $h$ so that%
\begin{align}
R_{N}  &  \leq C\{\sum_{\substack{K+L=N,\\K,L\leq N-1}}R_{K}R_{L}+\sum
_{l=2,p}[|\mathbf{w}^{N,h}|_{2,l}+R_{N-1}(\int_{r}^{T}|h(s)|^{l}%
ds)^{1/l}\label{eq7}\\
&  +1_{N=1}(\int_{r}^{T}|h(s)|^{2l}ds)^{1/2l}(\int_{r}^{T}||\mathbf{g}%
(s)||_{1,l}^{2l}ds)^{1/2l}]\},\nonumber
\end{align}
where $R_{N}=\sup_{r\leq s\leq T}\left[  |\mathbf{v}^{N,h}(t)|_{2,p}%
+|\mathbf{v}^{N,h}(t)|_{2,2}\right]  .$
\end{lemma}

\begin{proof}
If $M=0$, then (\ref{pf3}) is Navier-Stokes equation (\ref{pf1}) and, by Lemma
\ref{leme} there is a unique $\mathbb{H}_{p}^{2}\cap\mathbb{H}_{2}^{2}$-valued
continuous solution $\mathbf{v}^{0}=\mathbf{u}_{0}$. We proceed by induction.
Assume there are unique $\mathbf{v}^{M}=\mathbf{v}^{M,h}\in C([0,T],\mathbb{H}%
_{p}^{2}\cap\mathbb{H}_{2}^{2})$ solving the system for $M=0,\ldots,N-1$.
Consider the equation for $\mathbf{v}^{N}$:%
\begin{align}
\mathbf{v}^{N}\left(  t\right)   &  =\mathbf{w}^{N,h}+\int_{r}^{t}%
\mathcal{P}[\partial_{i}\left(  a^{\,ij}\left(  s\right)  \partial
_{j}\mathbf{v}^{N}\left(  s\right)  \right)  +b^{i}(s)\partial_{i}%
\mathbf{v}^{N}(s)\label{eq1}\\
&  -v^{0,k}\left(  s\right)  \partial_{k}\mathbf{v}^{N}(s)-v^{N,k}\left(
s\right)  \partial_{k}\mathbf{v}^{0}(s)+\mathbf{f}_{N}^{h}(s)]ds\nonumber
\end{align}
with
\begin{align*}
\mathbf{f}_{N}^{h}(s)  &  =-\sum_{\substack{K+L=N,\\K,L\leq N-1}%
}v^{K,k}\left(  s\right)  \partial_{k}\mathbf{v}^{L}(s)+1_{M=1}(\mathbf{g}%
\left(  s\right)  ,h(s))_{Y}\\
&  +\left(  \sigma^{i}(s),h(s)\right)  _{Y}\partial_{i}\mathbf{v}^{N-1}(s)
\end{align*}

Since
\[
\int_{r}^{T}|\mathbf{f}_{N}^{h}(s)|_{1,l}^{l}ds<\infty,l=2,p,
\]
the existence and uniqueness follows by Proposition \ref{aprop1}. Also, by
Proposition \ref{aprop1}, for each $N\geq1$ there is a constant $C$
independent of $h$ such that, denoting $R_{N}^{h}=\sup_{s\leq T}\left[
|\mathbf{v}^{M,h}(s)|_{2,p}+|\mathbf{v}^{M,h}(s)|_{2,2}\right]  $. Since\ for
$l=2,p,$%
\begin{align*}
|\mathbf{f}_{N}^{h}(s)|_{1,l}  &  \leq C[\sum_{\substack{K+L=N,\\K,L\leq
N-1}}|\mathbf{v}^{K}\left(  s\right)  |_{2,l}|\mathbf{v}^{L}(s)|_{2,p}%
+|h(s)|_{Y}|\mathbf{v}^{M-1}(s)|_{2,l}\\
&  +|h(s)|_{Y}||\mathbf{g}(s)||_{1,l},
\end{align*}
the inequality (\ref{eq7}) follows by Cauchy-Scwarz inequality.
\end{proof}

\subsubsection{Proof of Theorem \ref{main}}

{Fix }$T<T_{1}${ and choose a special CONS }$\mathbf{b}=\left\{
e_{n}\right\}  \in\mathcal{B}$ {such that for each }$h\in L_{\infty}%
([0,T],Y)${ there is a sequence }$z(n)\in\mathcal{Z}${ (see Lemma
\ref{re1 copy(1)}) for which }$e_{z(n)}\rightarrow h${ in }$L_{p}\left(
[0,T],Y\right)  ,${ for all }$p\geq2${, as }$n\rightarrow\infty${ (for
example, taking in }$L_{2}\left(  [0,T]\right)  $ {a trigonometric basis or
unconditional }$L_{p}\left(  [0,T]\right)  ${-basis (Haar basis), see
\cite{ym}). By Lemma \ref{les1} with }$\xi_{\alpha}=\xi_{\alpha}(\mathbf{b}%
)${,}%
\[
\mathbf{w}|_{\mathcal{D}(\mathbf{b)}}=\sum_{\alpha}\mathbf{w}_{\alpha}%
\xi_{\alpha}.
\]
{According to Lemma \ref{l1}, there is a unique }$\mathcal{D}${-}%
$\mathbb{H}_{p}^{2}\cap\mathbb{H}_{2}^{2}${-solution }%
\[
\mathbf{u}(t,x)=\sum_{\alpha}\mathbf{u}_{\alpha}(t,x)\xi_{\alpha}\in
{\ C}\mathcal{D}^{\prime}\left(  \mathbf{b};[r,T],\mathbb{H}_{p}^{2}%
\cap\mathbb{H}_{2}^{2}\right)
\]
{of (\ref{nse0}) in }$[0,T].${ The coefficients }$\mathbf{u}_{\alpha}(t,x)${
satisfy (\ref{pr1}) and, by Remark \ref{c1}, (\ref{pf2}) holds for all }%
$M\geq0,z\in\mathcal{Z}${. Fix }$h\in L_{\infty}([0,T],Y)${ and consider an
arbitrary }$e_{z(n)}\rightarrow h${ in }$L_{p}\left(  [0,T],Y\right)  ,${ for
all }$p\geq2${, as }$n\rightarrow\infty${.}

Then
\begin{align*}
\mathbf{w}^{M,z(n)}  &  =\left\langle \mathbf{w},p_{M}(z(n))\right\rangle \\
&  =\sum_{|\alpha|=M}\mathbf{w}_{\alpha}(t,x)z(n)^{\alpha}/\sqrt{\alpha!}%
\in\mathbb{H}_{p}^{2}\cap\mathbb{H}_{2}^{2},
\end{align*}

{and }%
\begin{align*}
\mathbf{u}^{M,z(n)}(t,x)  &  =\left\langle \mathbf{u}(t,x),p_{M}%
(z(n))\right\rangle \\
&  =\sum_{|\alpha|=M}\mathbf{u}_{\alpha}(t,x)z(n)^{\alpha}/\sqrt{\alpha!}\in
C\left(  [0,T],\mathbb{H}_{p}^{2}\cap\mathbb{H}_{2}^{2}\right)  ,
\end{align*}
$M\geq0,$ is the unique solution to the system (\ref{pf3}) corresponding to
$h=e_{z(n)}$ and $\mathbf{w}^{M,h}=\mathbf{w}^{M,z(n)}.${ Recall, }%
$\mathbf{u}^{-1,z(n)}=\mathbf{0}$ and $\mathbf{u}^{0,z(n)}$ coincides with the
solution of Navier-Stokes equation $\mathbf{u}_{0}$ in Lemma \ref{leme}. {By
Lemma \ref{a1}, there is a unique }$\mathbf{v}^{M,h}\in C\left(
[r,T],\mathbb{H}_{p}^{2}\cap\mathbb{H}_{2}^{2}\right)  ,M\geq0$, solving
(\ref{pf3} with $\mathbf{w}^{M,h}=\left\langle \mathbf{w},p_{M}%
(h)\right\rangle \in\mathbb{H}_{p}^{2}\cap\mathbb{H}_{p}^{2}.$ We have
$\mathbf{v}^{-1,h}(t)=0,$ and $\mathbf{v}^{0,h}$ coincides with the solution
of Navier-Stokes equation $\mathbf{u}_{0}$ in Lemma \ref{leme}. By Lemma
\ref{a1} (see (\ref{eq7})), for every $M\geq1$ there is a constant $C_{0}$
independent of $n$ such that%
\begin{equation}
\sup_{r\leq t\leq T}\sum_{l=2,p}[|\mathbf{u}^{M,z(n)}(t)|_{2,l}\mathbf{+|v}%
^{M,h}(t)|_{2,l}]\leq C_{0}. \label{eq9}%
\end{equation}

For $\mathbf{V}_{n}^{M}=\mathbf{v}^{M,h}-\mathbf{u}^{M,z(n)},M\geq1,$ the
following equation holds ($M\geq1):$%
\begin{align*}
\mathbf{V}_{n}^{M}\left(  t\right)   &  =\mathbf{w}^{M}-\mathbf{w}%
^{M,z(n)}+\int_{r}^{t}\mathcal{P}[\partial_{i}\left(  a^{\,ij}\left(
s\right)  \partial_{j}\mathbf{V}_{n}^{M}\left(  s\right)  \right) \\
&  +[b^{i}(s)-u_{0}^{i}\left(  s\right)  ]\partial_{i}\mathbf{V}_{n}%
^{M}(s)+V_{n}^{i}(s)\partial_{i}\mathbf{u}_{0}(s)+\mathbf{G}_{n}(s)]ds,
\end{align*}
where $\mathbf{G}_{n}(s)=\mathbf{G}_{n}^{1}(s)+\mathbf{G}_{n}^{2}(s)$ with%
\[
\mathbf{G}_{n}^{1}(s)=-\sum_{\substack{K+L=M,\\K,L\geq1}}[V_{n}^{K,k}%
(s)\partial_{k}\mathbf{u}^{L,z(n)}(s)+u^{K,z(n),k}\left(  s\right)
\partial_{k}\mathbf{V}_{n}^{L}(s)]
\]
and%
\begin{align*}
\mathbf{G}_{n}^{2}(h,s)  &  =1_{M=1}(\mathbf{g}(s),h_{n}(s))_{Y}+\left(
\sigma^{i}(s),h(s)\right)  _{Y}\partial_{i}\mathbf{V}_{n}^{M-1,h}(s)\\
&  +\left(  \sigma^{i}(s),h_{n}(s)\right)  _{Y}\partial_{i}\mathbf{u}%
^{M-1,z(n)}(s).
\end{align*}

By Proposition \ref{aprop1} in Appendix II, for $L_{M}^{n}=\sup_{s\leq
T}\left[  |\mathbf{V}_{n}^{M}(t)|_{2,p}\mathbf{+|V}_{n}^{M}(t)|_{2,2}\right]
,$ with $M\geq1,$%
\[
L_{M}^{n}\leq C[A_{n}+\sum_{l=2,p}(\int_{r}^{T}|\mathbf{G}_{n}(s)|_{1,l}%
^{l}ds)^{1/l}],
\]
where $A_{n}=\sum_{l=2,p}|\mathbf{w}^{M}-\mathbf{w}^{M,z(n)}|_{2,l}.$ We
estimate
\begin{align*}
|\mathbf{G}_{n}^{1}(s)|_{1,l}  &  \leq CC_{0}\sum_{1\leq K\leq M-1}L_{K}%
^{n},\\
|\mathbf{G}_{n}^{2}(h,s)|_{1,l}  &  \leq C[1_{M=1}|h_{n}(s)|_{Y}%
|\mathbf{g}(s)|_{1,l}\\
&  +L_{M-1}^{n}+C_{0}|h_{n}(s)|],
\end{align*}
where $h_{n}=h-e_{z(n)}$. So, for each $M\geq1$ there is a constant
independent of $n$ such that%
\begin{align*}
L_{M}^{n}  &  \leq C\{A_{n}+\sum_{1\leq K\leq M-1}L_{K}^{n}+\sum_{l=2,p}%
[\int_{0}^{T}|h_{n}(s)|^{l}ds)^{1/l}\\
&  +(\int_{0}^{T}|h_{n}(s)|^{2l}ds)^{1/2l}(\int|\mathbf{g}(s)|_{1,l}%
^{2l}ds)^{1/2l}]\}.
\end{align*}
Starting with $M=0,L_{0}^{n}=0$ for all $n$, it follows by induction that%
\[
L_{M}^{n}=\sup_{r\leq s\leq T}\sum_{l=2,p}|\mathbf{v}^{M,h}(s)-\mathbf{u}%
^{M,z(n)}(s)|_{2,l}\rightarrow0
\]
as $n\rightarrow\infty,M\geq1$.

{Since }$h\in L_{\infty}([0,T],Y)${ is arbitrary, }%
\[
\mathbf{u}(t)=\sum_{\alpha}\mathbf{u}_{\alpha}(t)\xi_{\alpha}\in
C\mathcal{T}^{\prime}(\mathbf{b};[r,T],\mathbb{H}_{p}^{2}\cap\mathbb{H}%
_{2}^{2})\subseteq C\mathcal{T}^{\prime}([r,T],\mathbb{H}_{p}^{2}%
\cap\mathbb{H}_{2}^{2})
\]
{ is the unique solution to (\ref{nse}). Obviously,}%
\[
\mathbf{v}^{M,h}(t,x)=\left\langle p_{M,T}(h),\mathbf{u}(t,x)\right\rangle
,h\in L_{\infty}([0,T],Y)
\]
satisfies (\ref{pf3}). Since for each CONS $\mathbf{b}^{\prime}=(e_{k}%
^{\prime})\in\mathcal{B}${ any linear combination of }$e_{k}^{\prime}${
belongs to }$L_{\infty}([0,T],Y),${ the generalized solution }$\mathbf{u}(t)${
extends any }$\mathcal{D}${-solution.}

{Now we will prove that the unique generalized }$\mathbb{H}_{p}^{2}%
\cap\mathbb{H}_{2}^{2}${-solution of (\ref{nse}) in }$[0,T]${ is }%
$\mathbb{F}^{W}${-adapted. We fix }$t^{\ast}\in(0,T),r\leq t^{\ast}${ and
consider a special basis }$\mathbf{\bar{b}}\in\mathcal{B}$ with $\bar{m}%
_{i}(t)${ in }$L_{2}((0,T))${ so that each }$\bar{m}_{i}${ is supported either
in }$[0,t^{\ast}]${ or in }$[t^{\ast},T]${ and such that for each }$h\in
L_{\infty}([0,T],Y)${ there is a sequence }$z(N)\in\mathcal{Z}${ (see Lemma
\ref{re1 copy(1)}) for which }$e_{z(N)}\rightarrow h${ in }$L_{p}\left(
[0,T],Y\right)  ,${ for all }$p\geq2${, as }$N\rightarrow\infty${ (for
example, }$\left(  \bar{m}_{k}\right)  ${ is a combination of two
trigonometric or unconditional }$L_{p}\left(  [0,T]\right)  ${-basis (Haar
basis) on }$(0,t^{\ast})${ and }$(t^{\ast},T)${). Let }$\bar{\xi}_{\alpha
}=\bar{\xi}_{\alpha}(\mathbf{\bar{b}}),\ \alpha\in\mathcal{J}${, the
corresponding orthonormal basis in }$L_{2}(\mathbb{F}_{T}^{W})${. Let }%
\[
\mathbf{u}(t)=\sum_{\alpha}\mathbf{\bar{u}}_{\alpha}(t)\bar{\xi}_{\alpha}\in
C\mathcal{T}^{\prime}(\mathbf{\bar{b}};[0,T],\mathbb{H}_{p}^{2}\cap
\mathbb{H}_{2}^{2})
\]
be the unique solution to (\ref{nse}) constructed using the representation
\[
\mathbf{w|}_{\mathcal{D}(\mathbf{\bar{b}})}=\sum_{\alpha}\mathbf{\bar{w}%
}_{\alpha}\bar{\xi}_{\alpha}\in\mathcal{T}^{\prime}(\mathbf{\bar{b}%
};\mathbb{H}_{p}^{2}\cap\mathbb{H}_{2}^{2}).
\]
So, $\mathbf{\bar{u}}_{\alpha}\in C\left(  [0,T],\mathbb{H}_{p}^{2}%
\cap\mathbb{H}_{2}^{2}\right)  $ satisfy (\ref{pr1}) in $\mathbb{L}%
_{l},l=2,p,$ with $\mathbf{w}_{\alpha}=\mathbf{\bar{w}}_{\alpha},\alpha\in J$.
{Let }$J^{\prime}=${\{}$\alpha\in J:\alpha${ has a non zero component
corresponding to }$\bar{m}_{k}${ whose support is in }$(t^{\ast},T)${\}. Since
}${\mathbf{w}}$ is $\mathcal{F}_{r}^{W}$-measurable, $\mathbf{\bar{w}}%
_{\alpha}=0$ if $\alpha\in J^{\prime}$. Indeed, if $\alpha\in J^{\prime}$,
there are $\ $\ $c_{i}\in\mathbf{R,}$ $z_{i}\in\mathcal{Z},i,=1\ldots,n,$ so
that
\[
\xi_{\alpha}=\sum_{i=1}^{n}c_{i}p_{T}(z_{i}).
\]
Then%
\begin{align*}
\mathbf{\bar{w}}_{\alpha}  &  =\left\langle \mathbf{w},\sum_{i=1}^{n}%
c_{i}p_{T}(z_{i})\right\rangle =\sum_{i=1}^{n}c_{i}\left\langle \mathbf{w}%
,p_{T}(z_{i})\right\rangle \\
&  =\sum_{i=1}^{n}c_{i}\left\langle \mathbf{w},p_{r}(z_{i})\right\rangle
=\left\langle \mathbf{w},\sum_{i=1}^{n}c_{i}p_{r}(z_{i})\right\rangle =0,
\end{align*}
because, by (\ref{Eq:WeAlpha}), { }%
\[
\sum_{i=1}^{n}c_{i}p_{r}(z_{i})=\mathbf{E}[\sum_{i=1}^{n}c_{i}p_{T}%
(z_{i})|\mathcal{F}_{r}^{W}]=\mathbf{E}\left[  \xi_{\alpha}|\mathcal{F}%
_{r}^{W}\right]  =0.
\]
We claim that, similarly, f{or }$t\in(0,t^{\ast})${, }$\mathbf{\bar{u}%
}_{\alpha}(t)=0${ if }$\alpha\in J^{\prime}${.} Indeed if $\alpha\in
J^{\prime}$ and $|\alpha|=1,$ we have $\mathbf{\bar{w}}_{\alpha=}0$ and in the
equation (\ref{pf4}) for $\mathbf{\bar{u}}_{\alpha}$ we have the input
function $\mathbf{F}_{\alpha}(t)=0$ for $t\in(0,t^{\ast})$. Therefore the
unique solution $\mathbf{\bar{u}}_{\alpha}(t)=0$ if $t\in(0,t^{\ast})$. Then
we simply apply induction on $|\alpha|=n$ and use (\ref{pf4}) (note that if
$|\alpha|=n+1,$ then $\alpha=\tilde{\alpha}+\varepsilon_{k}$ for some $k$ and
without any loss of generality we can assume that $\tilde{\alpha}\in
J^{\prime})$. {As a result,}%
\[
\mathbf{u}(t)=\sum_{\alpha\in\mathcal{J}}\mathbf{\bar{u}}_{\alpha}(t)\bar{\xi
}_{\alpha}=\sum_{\alpha\notin\mathcal{J}^{\prime}}\mathbf{\bar{u}}_{\alpha
}(t)\bar{\xi}_{\alpha},\ t\in\lbrack r,t^{\ast}].
\]
{Obviously, }$\bar{\xi}_{\alpha}${ are }$\mathcal{F}_{t^{\ast}}^{W}%
${-measurable for }$\alpha\notin J^{\prime}${. Also, for any }$z\in
\mathcal{Z},M\geq1,t\leq t^{\ast},$%
\begin{align*}
\left\langle p_{M}(z),\mathbf{u}(t)\right\rangle  &  =\sum_{\alpha
\notin\mathcal{J}^{\prime}}\mathbf{\bar{u}}_{\alpha}(t)\frac{z^{\alpha}}%
{\sqrt{\alpha!}}=\left\langle p_{M,t^{\ast}}(z),\mathbf{u}(t)\right\rangle \\
&  =\left\langle p_{M}(e_{z}1_{(0,t^{\ast})}),\mathbf{u}(t)\right\rangle
\end{align*}
{(note that }$e_{z}=\sum_{k}z_{k}e_{k},e_{z}1_{(0,t^{\ast})}=\sum_{k\notin
G}z_{k}e_{k}${, where }$G${ is the set of all }$k${ such that }$\bar{m}%
_{j_{k}}${ in }$e_{k}=\bar{m}_{j_{k}}l_{j_{k}}${ has its support in }%
$(t^{\ast},T))${. The statement of Theorem \ref{main} is proved.}

{The solution above has the restarting property as well. By the same arguments
as in Corollary \ref{cornew1} we have}

\begin{corollary}
{\label{cornew2}Let }$\mathbf{w}\in\mathcal{T}^{\prime}(\mathbb{H}_{2}^{2}%
\cap\mathbb{H}_{2}^{2})$ and {A1-A3 hold. Let }$u^{r,\mathbf{w}}(t)${ be
}$\mathbb{H}_{p}^{2}\cap\mathbb{H}_{2}^{2}${-solution to (\ref{nse})in
}$[r,T],T<T_{1},${ and }$r\leq r^{\prime}\leq t\leq T.${ Then}%
\begin{equation}
\mathbf{u}^{r,\mathbf{w}}(t)=\mathbf{u}^{r^{\prime},\mathbf{u(}r^{\prime}%
)}(t). \label{for1a}%
\end{equation}

\end{corollary}

\begin{corollary}
{\label{MP} (Markov Property) Assume that the assumptions of Corollary
\ref{cornew2} hold true and, in addition, }$\mathbf{w}${ is }$\mathcal{F}%
_{r}^{W}-${measurable, then }$u^{r,\mathbf{w}}(t)${ is }$(\mathcal{F}_{t}%
^{W})_{t\geq r}${-adapted. \ This together with (\ref{for1a}) can be
interpreted as Markov property.}
\end{corollary}

Let choose a uniformly bounded basis $\mathbf{b}=\left\{  e_{n}\right\}  $ (
$\sup_{k,s}|e_{k}(s)|_{Y}<\infty$) and rescale the solution in Theorem
\ref{main}
\[
\mathbf{u}|_{\mathcal{D}(\mathbf{b})}=\sum_{\alpha}\mathbf{u}_{\alpha}%
(t)\xi_{\alpha},
\]
using the second quantization $C_{\varepsilon}$ operator in (\ref{eq0}) in the
Fock space $\mathcal{H}=\sum_{n}\mathcal{H}_{n}$ ($\mathcal{H}_{n}%
=\mathbf{H}^{\hat{\otimes}n},$ see Appendix I, 5.1). Recall,%
\[
C_{\varepsilon}=\sum_{n=0}^{\infty}\kappa_{\varepsilon,n}D_{\varepsilon
}^{\otimes n},
\]
where $D_{\varepsilon}e_{k}=2^{-\varepsilon k}e_{k}${ and }$\kappa
_{\varepsilon,n}=e^{-\varepsilon e^{n}}$. {According to Proposition \ref{p1}},%
\begin{equation}
\mathbf{u}^{\varepsilon}(t)=C_{\varepsilon}\mathbf{u}(t)=\sum_{n=0}^{\infty
}\sum_{|\alpha|=n}\kappa_{\varepsilon,n}(2^{-\varepsilon\mathbf{N}})^{\alpha
}\mathbf{u}_{\alpha}(t)\xi_{\alpha} \label{eq-1}%
\end{equation}
is $\mathbb{H}_{2}^{2}$-valued continuous. The following statement holds.

\begin{proposition}
\label{p2}Let $\mathbf{b}=\left\{  e_{n}\right\}  $ be uniformly bounded and
A1-A3 hold. Let $\mathbf{u}\in\mathcal{T}^{\prime}([0,T],\mathbb{H}_{p}%
^{2}\cap\mathbb{H}_{2}^{2})$ be the solution to (\ref{nse0}) in $[0,T]$ with
the deterministic initial $\mathbf{w}\in\mathbb{H}_{p}^{2}\cap\mathbb{H}%
_{2}^{2}$. Then the rescaled $\mathbb{H}_{2}^{2}$-valued square integrable
continuous process $\mathbf{u}^{\varepsilon}(t)$ (defined by (\ref{eq-1}), see
Proposition \ref{p1}) is adapted and Markov (in a standard, rather than
generalized, sense).
\end{proposition}

\begin{proof}
For each $M\geq1,z\in\mathcal{Z},$%
\begin{align*}
\mathbf{E}p_{T,M}(z)\mathbf{u}^{\varepsilon}(t)  &  =\left\langle
p_{T,M}(z),\mathbf{u}^{\varepsilon}(t)\right\rangle =\kappa_{\varepsilon
,M}\left\langle p_{T,M}(z^{\varepsilon}),\mathbf{u}(t)\right\rangle \\
&  =\kappa_{\varepsilon,M}\left\langle p_{t,M}(z^{\varepsilon}),\mathbf{u}%
(t)\right\rangle =\left\langle p_{t,M}(z),\mathbf{u}^{\varepsilon
}(t)\right\rangle \\
&  =\mathbf{E}p_{t,M}(z)\mathbf{u}^{\varepsilon}(t).
\end{align*}
Since $\mathbf{u}^{\varepsilon}(t)$ is square integrable, $\mathbf{E}\left[
\mathbf{u}^{\varepsilon}(t)|\mathcal{F}_{t}^{W}\right]  =\mathbf{u}%
^{\varepsilon}(t)$. So $\mathbf{u}^{\varepsilon}(t)$ is adapted in a standard sense.

For any $0\leq s\leq t,$ by Corollary \ref{cornew2}, $\mathbf{u}%
(t)=\mathbf{u}^{s,\mathbf{u(}s)}(t)$. Therefore,%
\[
\mathbf{u}^{\varepsilon}(t)=C_{\varepsilon}\mathbf{u}(t)=C_{\varepsilon
}(\mathbf{u}^{s,\mathbf{u(}s)}(t))=(\mathbf{u}^{\varepsilon})^{s,\mathbf{u}%
^{\varepsilon}(s)}(t)
\]
and the standard Markov property follows.
\end{proof}

\begin{acknowledgement}
We are very grateful to S. Kaligotla and S. Lototsky for useful discussion.
\end{acknowledgement}

\section{{\label{Quant}}Appendix I. Wiener chaos}

In the first part of Appendix we present some facts of white noise analysis.

\subsection{{Rescaling of Wiener chaos by second quantization operator}}

{Consider a generalized random variable }$u=\sum_{\alpha}u_{\alpha}\xi
_{\alpha}=\sum_{\alpha}u_{\alpha}\sqrt{|\alpha|!}W(e_{\alpha})\in
\mathcal{D}^{\prime}(\mathbf{b}),\mathbf{b}=\left\{  e_{k},k\geq1\right\}  \in
B${, where }%
\[
W(e_{\alpha})=W^{\otimes n}(e_{\alpha})=\int_{0}^{T}\int_{0}^{s_{n}}\ldots
\int_{0}^{s_{2}}e_{\alpha}(s_{1},\ldots,s_{n})dW_{s_{1}}\ldots dW_{s_{n}},
\]
if $|\alpha|=n$ and $\left\{  e_{\alpha},\alpha\in J\right\}  $ defined by
(\ref{f0}){ is a CONS of the symmetric part }$\mathcal{H}_{n}=\mathbf{H}%
^{\hat{\otimes}n}${ of }$\mathbf{H}^{\otimes n}${ (recall }$\mathbf{H}%
=L_{2}(0,T])\times Y${)}$.$ {We can interpret }%
\[
u=\sum_{\alpha}u_{\alpha}\xi_{\alpha}=\sum_{n=0}^{\infty}\sum_{|\alpha
|=n}u_{\alpha}\sqrt{n!}W^{\otimes n}(e_{\alpha})
\]
{as a result of the noise }$W${ acting on an element of the Fock space:
}$W(\hat{u})=u$ with{ }%
\begin{align*}
\hat{u}  &  =\sum_{\alpha}u_{\alpha}\sqrt{|\alpha|!}e_{\alpha}=\sum
_{n=0}^{\infty}\sum_{|\alpha|=n}u_{\alpha}\sqrt{n!}e_{\alpha}\\
&  =\sum_{n=0}^{\infty}\hat{u}^{(n)}\in\mathcal{H}=\sum_{n=0}^{\infty
}\mathcal{H}_{n}=\sum_{n=0}^{\infty}\mathbf{H}^{\hat{\otimes}n}.
\end{align*}
{Here }$\mathbf{H}^{\hat{\otimes}0}=\mathbf{R}$ and the norm in the Fock space
$\mathcal{H}$ is defined as
\[
||\hat{u}||_{\mathcal{H}}^{2}=||\sum_{n=0}^{\infty}\hat{u}^{(n)}%
||_{\mathcal{H}}^{2}=\sum_{n=0}^{\infty}\frac{||\hat{u}^{(n)}||_{\mathcal{H}%
_{n}}^{2}}{n!}.
\]
Obviously, $\mathbf{E}[W(\hat{u})^{2}]=||\hat{u}||_{\mathcal{H}}^{2}.$ {Let
}$A=\left(  A_{n}\right)  _{n\geq0.}${ be a self-adjoint positive operator in
}$\mathcal{H}${ such that }$A_{n}e_{\alpha}=\lambda\left(  \alpha\right)
e_{\alpha}${, where }$\left\vert \alpha\right\vert =n${ and }$\lambda\left(
\alpha\right)  ${, }$\alpha\in J${, are positive numbers.}

\begin{remark}
The operator $A$ in the Fock space $\mathcal{H}$ can be used to rescale a
generalized r.v. For
\[
u=\sum_{\alpha}u_{\alpha}\xi_{\alpha}=W(\hat{u})\in\mathcal{D}^{\prime
}(\mathbf{b}),
\]
we define $Au=u^{A}\in\mathcal{D}^{\prime}(\mathbf{b)}$ by
\begin{align*}
Au  &  =u^{A}=W(A\hat{u})=\sum_{\alpha}u_{\alpha}\sqrt{|\alpha|!}W(Ae_{\alpha
})\\
&  =\sum_{\alpha}u_{\alpha}\lambda(\alpha)\sqrt{|\alpha|!}W(e_{\alpha}%
)=\sum_{\alpha}u_{\alpha}\lambda(\alpha)\xi_{\alpha}.
\end{align*}

\end{remark}

\begin{definition}
{Since }$\lambda(\alpha)>0${, we can define}%
\[
A^{-1}u=u^{A^{-1}}=\sum_{\alpha}u_{\alpha}\lambda(\alpha)^{-1}\xi_{\alpha}.
\]

\end{definition}

\begin{example}
{\label{ex3}1. (second quantization operator in space-time) Consider a
self-adjoint positive operator }$B${ in }$\mathbf{H}${ such that }%
$Be_{k}=\lambda_{k}e_{k}${. The second quantization operator }${A=}%
\Gamma(B)=\left(  B^{\otimes n}\right)  $ {in $\mathcal{H}$ is defined as }%
\[
Ae_{\alpha}=\Gamma(B)e_{\alpha}=B^{\otimes n}e_{\alpha}=\lambda^{\alpha
}e_{\alpha},|\alpha|=n,
\]
{where }$\lambda=\left(  \lambda_{k}\right)  ${ and }$\lambda^{\alpha}=\Pi
_{k}\lambda_{k}^{\alpha_{k}}${. We have}%
\[
\Gamma(B)u=\sum_{\alpha}u_{\alpha}\lambda^{\alpha}\xi_{\alpha}.
\]

{2. (second quantization operator in space) Consider a self-adjoint positive
operator }$B${ on }$Y${ such that the sequence of its eigenvectors }$\left(
\ell_{p}\right)  _{p\geq1}${ (}$B\ell_{p}=\lambda_{p}\ell_{p},\lambda_{p}>0)${
is a CONS in }$Y${. Let }$b=\left\{  e_{k},k\geq1\right\}  ,${ where }%
$e_{k}(s)=m_{i_{k}}(s)\ell_{j_{k}}.${ We extend }$B${ to }$\mathbf{H}${ by}%
\[
Be_{k}=B(m_{i_{k}}\ell_{j_{k}})=m_{i_{k}}B\ell_{j_{k}}=\lambda_{j_{k}}e_{k}%
\]
{and rescale in space-time using }$A=\left(  B^{\otimes n}\right)  ${. For
}$u=\sum_{\alpha}u_{\alpha}\xi_{\alpha}${ we have}%
\[
\Gamma(B)u=\sum_{\alpha}u_{\alpha}\lambda^{\alpha}\xi_{\alpha},
\]
{where }$\lambda^{\alpha}=\Pi_{k}\lambda_{j_{k}}^{\alpha_{k}}.$

{3. Consider a self-adjoint positive operator }$B${ on }$H${ such that
}$Be_{k}=\lambda_{k}e_{k}${ and a sequence of positive numbers }$q_{n}${. Let
}$A=\sum_{n=0}^{\infty}q_{n}B^{\otimes n}${. Then }$Ae_{\alpha}=q_{n}%
B^{\otimes n}e_{\alpha}=q_{n}\lambda^{\alpha}e_{\alpha}$, $|\alpha|=n.${In
this case,}%
\[
Au=\sum_{n=0}^{\infty}q_{n}\sum_{|\alpha|=n}u_{\alpha}\lambda^{\alpha}%
\xi_{\alpha}.
\]

\end{example}

{For the Wick product we have the following obvious statement.}

\begin{remark}
{\label{re4}Assume }$A=(A_{n})${ is a self-adjoint positive operator on
}$\mathcal{H}${ such that }$Ae_{\alpha}=A_{n}e_{\alpha}=\lambda(\alpha
)e_{\alpha},|\alpha|=n${ and }$\lambda(\alpha)${,}$\alpha\in I,${ are positive
numbers, }$u,v\in\mathcal{D}^{\prime}(\mathbf{b})${. Then, denoting }%
$Au=u^{A},Av=v^{A},${ we have}%
\[
A(u\lozenge v)=A\left(  A^{-1}u^{A}\lozenge A^{-1}v^{A}\right)  =\sum_{\alpha
}c_{\alpha}\xi_{\alpha},
\]
{where}%
\[
c_{\alpha}=\sum_{\beta\leq\alpha}\frac{\lambda(\alpha)}{\lambda(\beta
)\lambda(\alpha-\beta)}u_{\beta}^{A}v_{\alpha-\beta}^{A}\sqrt{\frac{\alpha
!}{\beta!(\alpha-\beta)!}}.
\]
{In particular, if \ }$\Gamma(B)${~is the second quantization operator in
space-time, then}%
\[
\frac{\lambda(\alpha)}{\lambda(\beta)\lambda(\alpha-\beta)}=\frac
{\lambda^{\alpha}}{\lambda^{\beta}\lambda^{\alpha-\beta}}=1
\]
{and }$\Gamma(B)(u\lozenge v)=u^{B}\lozenge v^{B}.$

{2. For the Skorokhod stochastic integral, we have }%
\begin{align*}
A\left(  \delta(u)\right)   &  =A\int_{0}^{T}u_{s}\lozenge\dot{W}_{s}%
ds=\int_{0}^{T}A\left(  A^{-1}u_{s}^{A}\lozenge A^{-1}\dot{W}_{s}^{A}\right)
ds\\
&  =\int_{0}^{T}A\left(  A^{-1}u_{s}^{A}\lozenge\dot{W}_{s}\right)  ds;
\end{align*}
{for the coefficients}%
\begin{align*}
\left(  A\delta(u)\right)  _{\alpha}  &  =\sum_{k}\int_{0}^{T}(u_{\alpha
(k)}^{A}(t),\lambda(\varepsilon_{k})e_{k}(t))_{Y}dt\sqrt{\alpha_{k}}%
\frac{\lambda(\alpha)}{\lambda(\alpha(k))\lambda(\varepsilon_{k})}\xi_{\alpha
}\\
&  =\sum_{k}\int_{0}^{T}(u_{\alpha(k)}^{A}(t),e_{k}(t))_{Y}dt\sqrt{\alpha_{k}%
}\frac{\lambda(\alpha)}{\lambda(\alpha(k))}\xi_{\alpha},\\
\left(  A(u(t)\lozenge\dot{W}_{t})\right)  _{\alpha}  &  =\sum_{k}%
(v_{\alpha(k)}^{A}(t),e_{k}(t))_{Y}\sqrt{\alpha_{k}}\frac{\lambda(\alpha
)}{\lambda(\alpha(k))}\xi_{\alpha}.
\end{align*}

\end{remark}

\subsection{\label{SSpwp}Product, Wick product, and Malliavin derivatives}

Consider Hilbert space $\mathbf{H}=L_{2}([0,T],Y),${ CONS }$\mathbf{b}%
=\left\{  e_{k},k\geq1\right\}  \in\mathcal{B},$ cylindrical Brownian motion
$W_{t}$, and Cameron-Martin basis $\left\{  \xi_{\alpha}\right\}  _{\alpha\in
J}$ introduced in Section \ref{ss1}. { }The \emph{Malliavin derivative}
$\mathcal{D}$ (see e.g. \cite{nualart}) is defined in $\mathcal{D}%
(\mathbf{b)}$ as follows (it assigns to $\xi_{\alpha}$ an element of
$\mathcal{D}\mathbf{(b};\mathbf{H)}$):
\begin{equation}
\mathcal{D}(\xi_{\boldsymbol{\mu}})=\sum_{k\geq1}\sqrt{\mu_{k}}\,\xi
_{\boldsymbol{\mu}-\boldsymbol{\varepsilon}_{k}}\,e_{k}=\sum_{\alpha}\sum
_{\mu=\alpha+\varepsilon_{k}}\sqrt{\mu_{k}}e_{k}\xi_{\alpha}. \label{eq:Dxi}%
\end{equation}
By induction,
\begin{equation}
\mathcal{D}^{n}(\xi_{\mu})=\sum_{\alpha}\left(  \sum_{|p|=n}1_{p+\alpha=\mu
}\sqrt{\frac{(\alpha+p)!}{\alpha!}}u_{p}\right)  \xi_{\alpha}, \label{Dnxialp}%
\end{equation}
where $u_{p}=\sum_{k_{1},\ldots,k_{n}}\sum_{\varepsilon_{k_{1}}+\ldots
+\varepsilon_{k_{_{n}}}=p}e_{k_{1}}\otimes\ldots\otimes e_{k_{n}}\in
\mathbf{H}^{\otimes n}.$

\begin{proposition}
\label{prop:wickvsdot}Let $\xi_{\theta}$ and $\xi_{\kappa}$ be elements of the
Cameron-Martin basis. Then, with probability 1,
\begin{equation}
\xi_{\theta}\xi_{\kappa}=\sum_{n=0}^{\infty}\frac{\mathcal{D}^{n}\xi_{\theta
}\diamondsuit\mathcal{D}^{n}\xi_{\kappa}}{n!}\text{ .} \label{eq:thetakap}%
\end{equation}

\end{proposition}

\begin{proof}
\bigskip It is a standard fact (see e.g. \cite{PAM}) that%
\[
\xi_{\theta}\xi_{\kappa}=\sum_{p\leq\theta\wedge\kappa}\sqrt{\binom{\theta}%
{p}\binom{k}{p}\binom{\theta+\kappa-2p}{\kappa-p}}p!\xi_{\theta+\kappa-2p}.
\]
Let us rewrite this expression as follows:
\[
\xi_{\theta}\xi_{\kappa}=\sum_{\substack{p,\beta,\gamma:\\p+\gamma
=\kappa,p+\beta=\theta}}\frac{\sqrt{\theta!\kappa!(\beta+\gamma)!}}%
{p!(\beta)!(\gamma)!}\xi_{\beta+\gamma}%
\]
where the summation goes over all triples $(p,\beta,\gamma)\in J\times J\times
J$ such that $p+\gamma=\kappa,p+\beta=\theta$. Changing variables (1-to-1
mapping that assigns to $(p,\beta,\gamma)$ the vector ($p,\beta,\alpha)$ with
$\alpha\geq\beta$) of summation by $p=p,\beta=\beta$, $\gamma+\beta=\alpha$,
we get%
\begin{align}
\xi_{\theta}\xi_{\kappa}  &  =\sum_{\substack{p,\,\beta\leq\alpha
:\\p+\alpha-\beta=\kappa,p+\beta=\theta}}\frac{\sqrt{\theta!\kappa!\alpha!}%
}{p!(\beta)!(\alpha-\beta)!}\xi_{\alpha}\label{eq: oldx}\\
&  =\sum_{\alpha}\sum_{\,\beta\leq\alpha}\sum_{n=0}^{\infty}\sum
_{|p|=n}1_{p+\alpha-\beta=\kappa}1_{p+\beta=\theta}\frac{\sqrt{\theta
!\kappa!\alpha!}}{p!(\beta)!(\alpha-\beta)!}\xi_{\alpha}\nonumber
\end{align}
By definition of the Wick product and (\ref{Dnxialp}), and taking into account
that $|p|!/p!$ is the number of different orthogonal unit vectors in $u_{p}%
\in$ $\mathbf{H}^{\otimes n}$, we arrive at%
\begin{align}
&  \mathcal{D}^{n}\xi_{\theta}\diamondsuit\mathcal{D}^{n}\xi_{\kappa
}\label{eq:newx}\\
&  =\sum_{\alpha}\sum_{\beta\leq\alpha}n!\sum_{|p|=n}1_{p+\alpha-\beta=\kappa
}1_{p+\beta=\theta}\frac{\sqrt{\theta!\kappa!\alpha!}}{p!(\beta)!(\alpha
-\beta)!}\xi_{\alpha}.\nonumber
\end{align}
Comparing (\ref{eq: oldx}) with (\ref{eq:newx}), we get (\ref{eq:thetakap}).
\end{proof}

\begin{remark}
\bigskip\label{re:hiorder} Proposition \ref{prop:wickvsdot} implies that
$\xi_{\theta}\xi_{\kappa}=\xi_{\theta}\diamondsuit\xi_{k}+\sum_{\gamma
<\theta+\kappa}c_{\gamma}\xi_{\gamma}.$ In other words, $\xi_{\theta
}\diamondsuit\xi_{k}=\xi_{\theta+k}$ is the \textit{highest stochastic order}
component of the Wiener chaos expansion of $\xi_{\theta}\xi_{\kappa}.$

By linearity, the statement of the Proposition could be extended to
\begin{equation}
XY=\sum_{n=0}^{\infty}\frac{(\mathcal{D}^{n}X)\lozenge(\mathcal{D}^{n}Y)}{n!},
\label{eq:pvswpfinite}%
\end{equation}
where $X$ and $Y$ are finite linear combinations of elements of Cameron-Martin
basis. If $X$ and $Y$ finite second moments, relation (\ref{eq:pvswpfinite})
could be derived from the former case by passing to the limit in $\ L_{1.}$By
linearity, the statement of Proposition could be extended to
\[
XY=\sum_{n=0}^{\infty}\frac{(\mathcal{D}^{n}X)\lozenge(\mathcal{D}^{n}Y)}%
{n!},
\]
where $X$ and $Y$ are finite linear combinations of elements of Cameron-Martin basis.
\end{remark}

\subsection{{Derivation of unbiased Navier-Stokes equation}}

{To simplify discussion, we will consider a velocity field which depends only
on one standard Gaussian random variable }$\eta\sim N\left(  0,1\right)  ,${
rather than a trajectory }o{f the Wiener process }$W_{t}${. }$\ ${An
interested reader would have little difficulties extending the arguments below
to the setting with Wiener process.}

{Consider a velocity field }%
\[
\mathbf{u}\left(  t,x\right)  =\sum_{n=0}^{\infty}\mathbf{u}_{n}\left(
t,x\right)  \xi_{n}\left(  \eta\right)  .
\]
{Note that in our setting the Cameron-Martin expansion (see Theorem
\ref{th:CM}) is indexed by integers rather than multi-indexes. Assume that for
every }$n${, }$u_{n}${ is analytic in }$x${ in that it could be written as}%
\[
\mathbf{u}_{n}\left(  t,x\right)  =\sum_{\gamma\in\mathbf{N}^{d}}%
\mathbf{c}_{n,\gamma}\left(  t\right)  x^{\gamma}.
\]
{Let }$Z=(Z_{1},\ldots,Z_{d})${ be a }$\mathcal{F}^{\eta}$-{measurable. Then
by substituting }$Z${ into }$u${ we get }%
\begin{equation}
\mathbf{u}\left(  t,Z\right)  :=\sum_{n}\left(  \sum_{\gamma}\mathbf{c}%
_{n,\gamma}\left(  t\right)  Z^{\gamma}\right)  \xi_{n}\left(  \eta\right)  .
\label{eq: yu}%
\end{equation}

{Now, let us introduce the Wick-powers of }$Z:$ $\ Z^{\lozenge\gamma}%
:=Z_{1}^{\diamondsuit\gamma_{1}}\lozenge...\lozenge Z_{d}^{\gamma_{d}}%
,\gamma=(\gamma_{1},\ldots,\gamma_{d})\in\mathbf{N}^{d}.$

{Next we will replace the standard algebra in (\ref{eq: yu}) by the Wick
algebra: \ }%
\[
\mathbf{u}_{n}^{\lozenge}\left(  t,Z\right)  :=\sum_{\gamma}\mathbf{c}%
_{n,\gamma}\left(  t\right)  Z^{\lozenge\gamma}%
\]
{Consider now the following random field}%
\[
\mathbf{u}^{\lozenge}\left(  t,Z\right)  :=\sum_{n\geq0}\mathbf{u}%
_{n}^{\lozenge}\left(  t,Z\right)  \lozenge\xi_{n}\left(  \eta\right)
\]

\begin{remark}
{\label{1a}Note that Wick algebra on nonrandom elements reduces to the
standard deterministic algebra.}
\end{remark}

{Let }$X_{t}=(X_{t}^{1},\ldots,X_{t}^{d})${ be a solution of the following
dynamic equation}%
\[
\dot{X}_{t}=\mathbf{u}^{\lozenge}\left(  t,X_{t}\right)  .
\]
{Then by the Wick chain rule}%
\begin{align*}
\overset{..}{X}_{t}  &  =\frac{d}{dt}\mathbf{u}^{\lozenge}\left(
t,X_{t}\right)  =\partial_{t}\mathbf{u}^{\lozenge}\left(  t,X_{t}\right)
+\nabla\mathbf{u}^{\lozenge}\left(  t,X_{t}\right)  \lozenge\dot{X}_{t},\\
&  =\partial_{t}\mathbf{u}^{\lozenge}\left(  t,X_{t}\right)  +\mathbf{u}%
^{\lozenge}\left(  t,X_{t}\right)  \nabla\lozenge\mathbf{u}^{\lozenge}\left(
t,X_{t}\right)  .
\end{align*}
{ If }$\mathbf{F}=\mathbf{F}(t,x)$ is an acting force, {this yields (Wick)
Euler equation }%
\[
\partial_{t}\mathbf{u}^{\lozenge}\left(  t,x\right)  =-\mathbf{u}^{\lozenge
}\left(  t,x\right)  \nabla\lozenge\mathbf{u}^{\lozenge}\left(  t,x\right)
+\mathbf{F}\left(  t,x\right)
\]
If there is no randomness, {due to Remark \ref{1a}, this equation reduces to
the standard Euler equation:}%

\[
\partial_{t}\mathbf{u}\left(  t,x\right)  =-\mathbf{u}\left(  t,x\right)
\nabla\mathbf{u}\left(  t,x\right)  +\mathbf{F}\left(  t,x\right)  .
\]
{Now, by taking }$\mathbf{F}=\Delta\mathbf{u}-\nabla P,${ where }$P${ stands
for pressure, we get the unbiased Navier-Stokes equation}%
\[
\partial_{t}\mathbf{u}^{\lozenge}\left(  t,x\right)  =\Delta\mathbf{u}%
-\mathbf{u}^{\lozenge}\left(  t,x\right)  \nabla\lozenge\mathbf{u}^{\lozenge
}\left(  t,x\right)  -\nabla P+\mathbf{F}\left(  t,x\right)  .
\]

\section{Appendix II. Stokes equation}

Consider a deterministic Stokes equation for $\mathbf{u=}(u^{l})_{1\leq l\leq
d},\,$\ and scalar functions $P,$
\begin{align*}
\partial_{t}\mathbf{u}\left(  t,x\right)   &  =\partial_{i}\left(
a^{\,ij}\left(  t,x\right)  \partial_{j}\mathbf{u}\left(  t,x\right)  \right)
+b^{i}(t,x)\partial_{i}\mathbf{u}(t,x)\\
+\mathbf{G(}t,x)\mathbf{u}(t,x)  &  +\mathbf{f}\left(  t,x\right)  +\nabla
P(t,x),\operatorname{div}\mathbf{u}(t)=0,\\
\mathbf{u}\left(  r,x\right)   &  =\mathbf{w}(x),x\in\mathbf{R}^{d},r\leq
t\leq T.
\end{align*}
equivalently,%
\begin{align}
\partial_{t}\mathbf{u}\left(  t\right)   &  =\mathcal{S}[\partial_{i}\left(
a^{\,ij}\left(  t\right)  \partial_{j}\mathbf{u}\left(  t\right)  \right)
+b^{i}(t)\partial_{i}\mathbf{u}(t)+\mathbf{G(}t)\mathbf{u}(t)+\mathbf{f}%
(t)\mathbf{]}\label{af1}\\
\mathbf{u}(r)  &  =\mathbf{w,}t\in\lbrack r,T],\nonumber
\end{align}
where $\mathcal{S}$ is the solenoidal projection of the vector fields,
\[
\mathbf{a}(t,x)=\left(  a^{ij}(t,x)\right)  _{1\leq i,j\leq d},\mathbf{b}%
(t,x)=\left(  b^{i}(t,x)\right)  _{1\leq i<d},\mathbf{G}(t,x)=\left(
g^{ij}(t,x)\right)  _{1\leq i,j\leq d}%
\]
are measurable bounded functions. The matrix $a$ is symmetric and positive.

We will need the following assumption.

\textbf{B}. For all $t\geq0,\,x,\lambda\in$\textbf{$R$}$^{d}$, $K|\lambda
|^{2}\geq a^{ij}(t,x)\lambda^{i}\lambda^{j}\geq\delta|\lambda|^{2},$where
$K,\delta$ are fixed strictly positive constants. Also, for all $(t,x)\in
\lbrack r,T]\times\mathbf{R}^{d},$%
\[
\max_{t,x,|\alpha|\leq2}|\partial_{x}^{\alpha}\mathbf{a}(t,x)|+\max
_{t,x,|\alpha|\leq1}|\partial_{x}^{\alpha}\mathbf{b}(t,x)|+\sup_{t,x}%
|\mathbf{G}(t,x)|\leq K.
\]

\begin{definition}
\label{def1 copy(1)}A function $\mathbf{u}\in C([r,T],\mathbb{H}_{p}%
^{s}\mathbb{(}\mathbf{R}^{d}))$ is an $\mathbb{H}_{p}^{s}$-solution of
(\ref{af1}) if the equality
\begin{equation}
\mathbf{u}(t)=\mathbf{w}+\int_{r}^{t}\mathcal{S}[\partial_{i}(a^{ij}%
(s)\partial_{j}\mathbf{u}(s)\mathbf{)+}b^{i}(s)\partial_{i}\mathbf{\mathbf{u}%
}(s)+\mathbf{f}(s)]ds, \label{sol1}%
\end{equation}
holds in $\mathbb{H}_{p}^{s-2}\mathbb{(}\mathbf{R}^{d})$ for every
$t\in\lbrack r,T].$
\end{definition}

\begin{proposition}
\label{aprop1}Let $p>d,$ assumption $\mathbf{B}$ hold, $\mathbf{w}%
\in\mathbb{H}_{p}^{2}\cap\mathbb{H}_{2}^{2},$%
\[
\int_{r}^{T}|\mathbf{f}(s)|_{1,l}^{l}ds+\int_{r}^{T}|\nabla\mathbf{G}%
(s)|_{l}^{l}ds<\infty,l=2,p.
\]
Then there is a unique $\mathbb{H}_{p}^{2}\cap\mathbb{H}_{2}^{2}$-valued
continuos solution to (\ref{af1}). Moreover, there is a constant $C$
independent of $\mathbf{f,w},\mathbf{u}$ so that%
\begin{align}
&  \sup_{r\leq t\leq T}[|\mathbf{u}(t)|_{2,2}+|\mathbf{u}(t)|_{2,p}%
]\label{afo4}\\
&  \leq C\left(  |\mathbf{w|}_{2,2}+|\mathbf{w|}_{2,p}+(\int_{r}%
^{T}|\mathbf{f}(s)|_{1,p}^{p}ds)^{\frac{1}{p}}+(\int_{r}^{T}|\mathbf{f}%
(s)|_{1,2}^{2}ds)^{1/2}\right)  .\nonumber
\end{align}

\end{proposition}

\begin{proof}
Let $\mathbf{w}\in\mathbb{H}_{p}^{3}\cap\mathbb{H}_{2}^{2}$. {By Proposition
4.7 and Corollary 4.6 in \cite{m} (applied for }$s=0)${, there is a unique
}$\mathbb{H}_{p}^{1}\cap\mathbb{H}_{2}^{1}$-valued continuous {solution
}$\mathbf{u}$ of (\ref{af1}), and%
\[
\int_{r}^{T}|\mathbf{u}(s)|_{2,l}^{l}ds<\infty,l=2,p\text{.}%
\]

Moreover,%
\begin{equation}
\sup_{s\leq T}|\mathbf{u}(s)|_{1,l}^{l}\leq C[|\mathbf{w|}_{2,l}^{l}+\int%
_{r}^{T}|\mathbf{f}(s)|_{l}^{l}ds],l=2,p, \label{eq2}%
\end{equation}

Consider Stokes equation%
\begin{align}
\mathbf{\xi}\left(  t\right)   &  =\mathbf{w}_{\alpha}+\int_{r}^{t}%
\mathcal{P}[\partial_{i}\left(  a^{\,ij}\left(  s\right)  \partial
_{j}\mathbf{\xi}\left(  s\right)  \right)  +\mathbf{F}(s)]ds\label{e2}\\
\text{ }\operatorname{div}\text{\thinspace}\mathbf{\xi}(t)  &  =0,t\in\lbrack
r,T],\nonumber
\end{align}
where%
\[
\mathbf{F}(s)=b^{i}(s)\partial_{i}\mathbf{u}(s)+\mathbf{G(}s)\mathbf{u}%
(s)+\mathbf{f}(s).
\]
It is readily checked (using Sobolev embedding theorem) that%
\begin{align}
|\mathbf{F}(s)|_{1,l}  &  \leq C[K(|\mathbf{u}(s)|_{2,l}+|\nabla
\mathbf{G}(s)|_{l}|\mathbf{u}(s)|_{1,p}\label{eq3}\\
+|\mathbf{f}(s)|_{1,l}],l  &  =2,p.\nonumber
\end{align}
and
\[
\int_{r}^{T}|\mathbf{F}(s)|_{1,l}^{l}ds<\infty,l=2,p.
\]
By Corollary \ 4.6 and Proposition 4.7 in \cite{m} (applied to (\ref{e2}) with
$s=1$), there is a unique $\mathbb{H}_{2}^{2}\cap\mathbb{H}_{p}^{2}$-valued
continuous solution of (\ref{e2}) $\mathbf{\xi=u}$ (by uniqueness) such that%
\[
\int_{r}^{T}|\mathbf{u}(s)|_{3,l}^{l}ds<\infty,l=2,p\text{.}%
\]

Let $\alpha$ be a multiindex such that $|\alpha|\leq2$. Then $\mathbf{u}%
_{\alpha}=\partial^{\alpha}\mathbf{u}$ is $\mathbb{L}_{p}\cap\mathbb{L}_{2}%
$-valued continuous and satisfies the equation
\begin{align*}
\partial_{t}\mathbf{u}_{\alpha}(t)  &  =\mathcal{S\{}\partial^{\alpha
}[\partial_{i}(a^{ij}(t)\partial_{j}\mathbf{u}(t)\mathbf{)+F}(t)\mathbf{],}\\
\mathbf{u}_{\alpha}(0)  &  =\partial^{\alpha}\mathbf{w}.
\end{align*}
Differentiating the product, we obtain
\begin{equation}
\partial^{\alpha}\partial_{i}(a^{ij}(t)\partial_{j}\mathbf{u}(t)\mathbf{)}%
=\partial_{i}(a^{ij}(t)\partial_{j}\mathbf{u}_{\alpha}(t)\mathbf{)+}%
\partial_{i}\mathbf{D}_{\alpha}(t) \label{eq4}%
\end{equation}
with
\begin{equation}
|\mathbf{D}_{\alpha}(t)|_{l}\leq C|\mathbf{u(}t)|_{2,l},l=2,p. \label{eq5}%
\end{equation}
By Lemma 3 in \cite{mr24}, $y_{l,\alpha}(t)=|\mathbf{u}_{\alpha}(t)|_{l}%
^{l},l=2,p,$ is differentiable:
\[
y_{l,\alpha}(t)=y_{l,\alpha}(r)+\int_{r}^{t}h_{l,\alpha}(s)\,ds,
\]
with
\begin{align*}
h_{l,\alpha}(s)  &  =l\{\langle|\mathbf{u}_{\alpha}(s)|^{l-2}\mathbf{u}%
_{\alpha}(s),\partial^{\alpha}\mathbf{F}(s)]\rangle_{1,l}\\
&  -\int a^{ij}(s)\partial_{i}(|\mathbf{u}_{\alpha}(s)|^{l-2}u_{\alpha}%
^{k}(s))\partial_{j}u_{\alpha}^{k}(s)\,dx\\
&  -\int\partial_{i}(|\mathbf{u}_{\alpha}(s)|^{l-2}u_{\alpha}^{k}%
(s))D_{\alpha}^{k}(s)\,dx\}.\,
\end{align*}
Notice $\partial^{\alpha}\mathbf{F}(s)\in\mathbb{H}_{-1,l}$ and, by our
assumptions, there is a constant $C$ so that for all $s\in\lbrack r,T]$
\begin{equation}
|\partial^{\alpha}\mathbf{F}(s)|_{-1,l}\leq C|\mathbf{F}(s)|_{1,l}%
,l=2,p.\nonumber
\end{equation}
We have $h_{l,\alpha}(s)=lh_{l,\alpha}^{1}(s)+lh_{l,\alpha}^{2}(s),$ where
\[
h_{l,\alpha}^{1}(s)=-\int a^{ij}(s)\partial_{i}(|\mathbf{u}_{\alpha}%
(s)|^{l-2}u_{\alpha}^{k}(s))\partial_{j}u_{\alpha}^{k}(s)\,dx.
\]
Then
\[
h_{l,\alpha}^{1}(s)\leq-\delta\int|\mathbf{u}_{\alpha}(s)|^{l-2}%
|\nabla\mathbf{u}(s)|^{2}\,dx,
\]
and for each $\varepsilon>0$ there is a constant $C_{\varepsilon}$ such that
\begin{align*}
|h_{l,\alpha}^{2}(s)|  &  \leq\varepsilon\int|\mathbf{u}_{\alpha}%
(s)|^{l-2}|\nabla\mathbf{u}_{\alpha}(s)|^{2}\,dx\\
&  +C_{\varepsilon}\int[|\mathbf{u}_{\alpha}(s)|^{l-2}(|\nabla\mathbf{F}%
(s)|^{2}+|\mathbf{D}_{\alpha}(s)|^{2})\\
&  +|\mathbf{u}_{\alpha}(s)|^{l-1}|\mathbf{F}(s)|]dx,
\end{align*}
So, we obtain that
\[
y_{l}(t)=\sum_{|\alpha|\leq2}y_{l,\alpha}(t)=y_{l}(r)+\int_{r}^{t}%
h_{l}(s)\,ds
\]

with
\begin{equation}
h_{l}(s)=\sum_{|\alpha|\leq2}h_{l,\alpha}(s)\leq C(y_{l}(s)+f_{l}%
(s)),\nonumber
\end{equation}
where $f_{l}(s)=|\mathbf{F(}s)|_{1,l}^{l}+|\mathbf{D}_{\alpha}(s)|_{l}^{l}.$
Therefore, by (\ref{eq5}), (\ref{eq3}),%
\[
h_{l}(s)\leq C[y_{l}(s)+|\mathbf{f}(s)|_{1,l}^{l}+|\nabla\mathbf{G}%
(s)|_{l}^{l}|\mathbf{u}(s)|_{1,p}^{l}.
\]
By Gronwall's inequality,
\[
\sup_{r\leq s\leq T}y_{p}(s)\leq C\mathbf{[}y_{p}(r)+\int_{r}^{T}%
|\mathbf{f}(s)|_{1,p}^{p}ds],
\]
and%
\begin{equation}
\sup_{r\leq t\leq T}|\mathbf{u}(t)|_{2,p}^{p}\leq C\left(  |\mathbf{w|}%
_{2,p}^{p}+\int_{r}^{T}|\mathbf{f}(s)|_{1,p}^{p}ds\right)  , \label{afo2}%
\end{equation}
where $C$ is independent of $\mathbf{w}$ and $\mathbf{f}$. Similarly, by
Gronwall's inequality%
\[
\sup_{r\leq s\leq T}y_{2}(s)\leq C\mathbf{[}y_{2}(r)+\sup_{s\leq T}%
|\mathbf{u}(s)|_{1,p}^{2}\int_{r}^{T}|\nabla\mathbf{G}(s)|_{2}^{2}ds+\int%
_{r}^{T}|\mathbf{f}(s)|_{1,2}^{2}ds],
\]
and ( see (\ref{eq2}))%
\begin{equation}
\sup_{r\leq t\leq T}|\mathbf{u}(t)|_{2,2}^{2}\leq C\left(  |\mathbf{w|}%
_{2,2}^{2}+|\mathbf{w|}_{2,p}^{2}+(\int_{r}^{T}|\mathbf{f}(s)|_{p}%
^{p}ds)^{\frac{2}{p}}+\int_{r}^{T}|\mathbf{f}(s)|_{1,2}^{2}ds\right)  ,
\label{afo3}%
\end{equation}
where $C$ is independent of $\mathbf{w},\mathbf{f}$ and $\mathbf{u}$.
Combining (\ref{afo2}) and (\ref{afo3}) we have (\ref{afo4}) with $C$ is
independent of $\mathbf{w},\mathbf{f}$ and $\mathbf{u}$.

Given $\mathbf{w}\in\mathbb{H}_{p}^{2}\cap\mathbb{H}_{2}^{2},$ there is a
sequence $\mathbf{w}_{n}\in\mathbb{H}_{p}^{3}\cap\mathbb{H}_{2}^{2}$ so that
$\mathbf{w}_{n}\rightarrow\mathbf{w}$ in $\mathbb{H}_{p}^{2}\cap\mathbb{H}%
_{2}^{2}$. For every $n$ there is a unique $\mathbb{H}_{p}^{2}\cap
\mathbb{H}_{2}^{2}$-valued continuous solution $\mathbf{u}_{n}$\ of
(\ref{sol1}) with the initial condition $\mathbf{u}(r)=\mathbf{w}_{n}$. By
(\ref{afo4})%
\[
\sum_{l=2,p}\sup_{r\leq t\leq T}|\mathbf{u}_{n}(t)-\mathbf{u}_{m}%
(t)|_{2,l}\leq C\sum_{l=2,p}|\mathbf{w}_{n}-\mathbf{w}_{m}|_{2,l}\rightarrow0
\]
as $n,m\rightarrow\infty.$ There is a continuous $\mathbb{H}_{p}^{2}%
\cap\mathbb{H}_{2}^{2}$-valued $\mathbf{u}(t)$ such that
\[
\sum_{l=2,p}\sup_{r\leq t\leq T}|\mathbf{u}_{n}(t)-\mathbf{u}(t)|_{2,l}%
\rightarrow0
\]
as $n\rightarrow\infty.$ Obviously, $\mathbf{u}$ is $\mathbb{H}_{p}^{2}%
\cap\mathbb{H}_{2}^{2}$-valued continuous solution of (\ref{sol1}) with
initial condition $\mathbf{u}(r)=\mathbf{w}$ and (\ref{afo4}) holds.{\ }
\end{proof}

\end{document}